\documentclass[11pt]{article}
\usepackage{amsmath}
\usepackage{amsfonts}
\usepackage{amssymb}
\usepackage{amsthm}
\usepackage{color}
\usepackage[greek,english]{babel}

\usepackage{graphicx}

\theoremstyle{definition}
\newtheorem{defn}{Definition}[section]
\theoremstyle{plain}

\newtheorem{lemma}[defn]{Lemma}
\newtheorem{tvr}[defn]{Proposition}
\newtheorem{cor}[defn]{Corollary}
\theoremstyle{remark}

\newcommand{\R}{\mathbb{R}}
\newcommand{\Z}{\mathbb{Z}}
\newcommand{\Q}{\mathbb{Q}}
\newcommand{\C}{\mathbb{C}}

\newcommand{\HH}{\mathbb{H}}
\newcommand{\F}{\mathbb{F}}
\newcommand{\Ztau}{\mathbb{Z}[\tau]}

\newcommand{\cS}{\mathcal{S}}
\newcommand{\cC}{\mathcal{C}}
\newcommand{\cD}{\mathcal{D}}
\newcommand{\cR}{\mathcal{R}}

\newcommand{\cU}{\mathcal{U}}

\newcommand{\card}{\mathrm{card}\,}
\allowdisplaybreaks
\usepackage{backref}

\begin{document}

\title{Discretization of SU(2) and the orthogonal group using icosahedral symmetries and the golden numbers}

\author{Robert~V.~Moody \& Jun Morita}

\date{\today}

\maketitle

\noindent
Department of Mathematics and Statistics, University of Victoria, Victoria, BC., V8W\,3R4 Canada
Department of Mathematics, University of Tsukuba, Tsukuba, Ibaraki, 305-8571, Japan.

\noindent $\phantom{^\S}$~E-mail: rvmoody@mac.com, morita@math.tsukuba.ac.jp

\bigskip
\begin{abstract} The vertices of the four-dimensional $600$-cell 
form a non-crystallographic root system whose corresponding symmetry group is the Coxeter group $H_{4}$. There are two special coordinate representations of this root system in which they and  their corresponding Coxeter groups involve only rational numbers and the golden ratio $\tau$. The two are related by the conjugation  $\tau \mapsto\tau' = -1/\tau$. This paper investigates what happens when the two root systems are combined
and the group generated by both versions of $H_{4}$ is allowed to operate on them. 
The result is a new, but infinite, `root system' $\Sigma$ which itself turns out to have a natural structure of the unitary group $SU(2,\cR)$ over the ring $\cR = \Z[\frac{1}{2},\tau]$ (called here golden numbers). Acting upon it is the naturally associated infinite reflection group $H^{\infty}$, which we prove is of index $2$ in the orthogonal group $O(4,\cR)$.
The paper makes extensive use of the quaternions over $\cR$ and leads to a highly structured discretized filtration of $SU(2)$. We use this to offer a simple and effective way to approximate any element of $SU(2)$ to any degree of accuracy required using the repeated actions of just five fixed reflections, a process that may find application in computational methods in quantum mechanics.   
\end{abstract}\


\smallskip
\section{Introduction}

The symmetry group of the icosahedron (and the dodecahedron) is the icosahedral group, denoted $H_{3}$ here, 
with $120$ elements. It is a finite Coxeter group, that is to say, it is a finite group generated by reflections and Coxeter 
relations,
and it is simply transitive on the simplicial cells into which the 
icosahedron is partitioned by the mirrors of the reflections in $H_3$. Apart from numerous dihedral groups, there are only two finite indecomposable non-crystallographic
Coxeter groups: $H_{3}$ and its big sister $H_{4}$, which is the symmetry group of the four-dimensional regular polytope called the $600$-cell (and its dual, the $120$-cell). This paper involves both these groups, but notably the latter. 

\begin{figure}\centering
\includegraphics[scale=0.35]{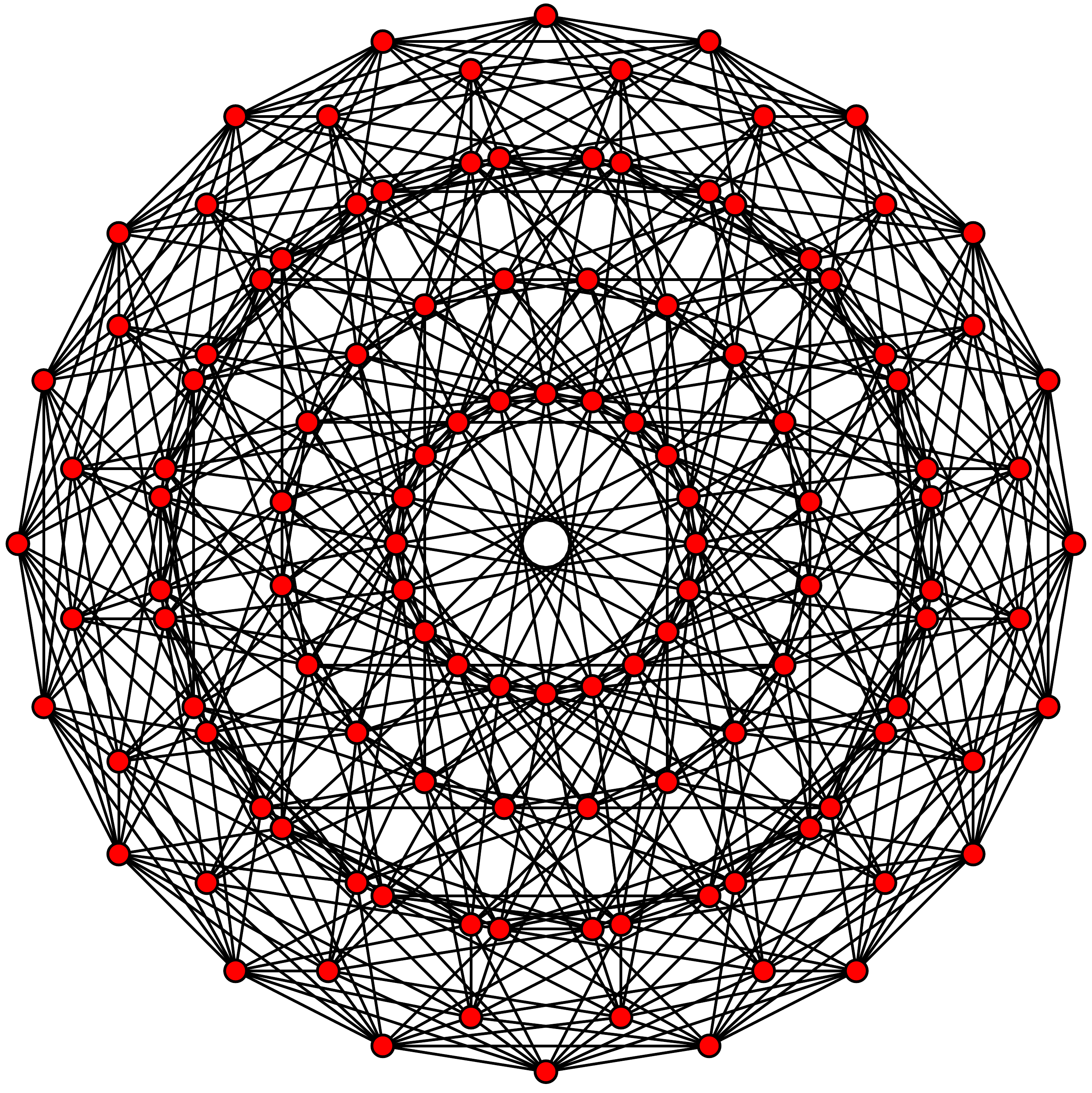}
\caption{A two-dimensional projection of the $600$-cell, whose vertices form a root system of type $H_{4}$. The dots make up the $120$ vertices that constitute a root system of type $H_4$.}
\label{600-cell}
\end{figure}

The $600$-cell has $120$ vertices and $600$ three-dimensional faces, see Fig.~\ref{600-cell}. The $120$ vertices form a root system, 
of type $H_{4}$. With this interpretation of the vertices as roots, the set of $60$ reflections
in these roots (opposite roots give the same reflection) is the entire set of reflections in the group
$H_{4}$, and they generate it. The order of the group is $|H_{4}| = 120^{2}=14400$. Notice that we use
the symbol $H_4$ both as an adjective signifying the type of root system involved and as a noun signifying
the reflection group generated by the roots of that system. We will do the same thing with $H_3$.

Both $H_{3}$ and $H_{4}$ involve the golden ratio
$\tau = (1  +\sqrt 5)/2$ and its algebraic conjugate $\tau' = (1  -\sqrt 5)/2$ in all sorts of significant ways. 
For instance, there is the well-known model of the vertices of the $600$-cell \cite{DV}\,\S20 and \S27, or equivalently the roots of
the root system of type $H_{4}$, as the set of points
\begin{equation}\label{120-roots}
\Delta =
\begin{cases}
\mbox{ {\bf R[0,0]} \quad $(\pm 1,0,0,0)$ and all its permutations: $8$ roots in all;}\\
\mbox{ {\bf R[1,0]}\quad  $\frac{1}{2}(\pm 1, \pm1, \pm1,\pm1)$: $16$ roots in all;}\\
\mbox{ {\bf R[1,1]} \quad $\frac{1}{2}(0, \pm1, \pm \tau', \pm \tau)$ and all its even permutations: $96$ roots in all}.
\end{cases}
\end{equation}

Here the first component of the label refers to the power of $2$ appearing in the denominators of the components
and the second component is used to distinguish the rational versus irrational nature of (some of) the components. 

One cannot help noticing the interesting fact that only even permutations are allowed in the third
type of root. Half of the potential permutations are missing. The other half can be obtained
by conjugating these roots by interchanging $\tau$ and $\tau'$ throughout. Of course doing so produces
another model $\Delta'$ of the root system of type $H_{4}$, and the reflections in these generate a second
model $H_{4}'$ of the Coxeter group $H_{4}$. We shall also use the notation 
$H_4= W(\Delta)$ and $H_4' = W(\Delta')$, that suggests
the origin of these groups as reflection groups.

It is tempting to look at the group generated by the reflections of
$\Delta \cup \Delta'$, but one quickly realizes that this group, let's call it $H^{\infty}$ (or $H^{\infty}[4]$) is
infinite.\footnote{Throughout the paper we will see parallel structures pertaining to the three-dimensional setting around $H_{3}$
and the four-dimensional setting around $H_{4}$. Generally we use the same symbol for the corresponding objects, but
if the context does not make it clear then we shall make the distinction as in $H^{\infty}[3]$, $H^{\infty}[4]$, 
$\Sigma[3]$, $\Sigma[4]$, $\cS[3]$, $\cS[4]$, etc.}
Since all the points
of the set
\begin{equation}\label{defSigma}
\Sigma= \Sigma[4]:=H^{\infty}(\Delta \cup \Delta')
\end{equation}
lie on the sphere of radius $1$ in $\R^{4}$, which is compact,
$\Sigma$ is certainly not a discrete set. Up to now  
no one seems to have paid much attention to it. 

The objective of this paper is to get some idea of what $\Sigma$ and $H^{\infty}$ look like. In fact they have some very attractive features, 
as we shall see. Not surprisingly $\Sigma$ is a dense set of points on the unit $3$-sphere $\cS$ in $\R^{4}$,
see Prop.~\ref{denseness}. 
It is also a group under quaternion multiplication, and we shall show that it is an amalgamation of
$\Delta$ and $\Delta'$, each of which is itself a group. Similarly $H^{\infty}$ is an amalgamation of
the two groups $H_4$ and $H_4'$. What is especially interesting is that using $\Sigma$ and  $H^{\infty}$ allows us to explicitly approximate elements of the two groups $SU(2)$ and $SO(3)$, the unitary and special orthogonal
groups of $\C^{2}$ and $\R^{3}$ respectively, by elements with matrix coefficients of the form $\frac{1}{2^{n}}\Z[\tau]$ (we call them dyadic integers from the  field $\Q[\tau]$). In fact the matrices
used all arise from the $96$ reflections of $\Delta \cup \Delta'$ which are the ones involving $\tau$ and $\tau'$,
and even from just $5$ reflections if one gets down to the level of root system bases. The approximation can be made as fine as one wishes (with increasing powers of two in the denominators) and there is a simple and efficient algorithm for doing it. 

The key to all of this is to interpret $\R^{4}$ as the standard division ring of quaternions, and use the fact that
its unit sphere can be identified with $SU(2)$ and that it can be made to act as $SO(3)$ on the three-dimensional
subspace of pure imaginary quaternions (those which change their sign under quaternionic conjugation).

The paper is primarily concerned with the $H_4$ picture, but we also require information from the corresponding,
and essentially parallel situation based on the three-dimensional icosahedral root system $H_3$
and its corresponding conjugate system $H_3'$. The arguments are not identical but sufficiently parallel that we provide
only sketches of the arguments for the three-dimensional case.

The paper is organized around understanding the structure of the root system $\Sigma$ and the group $H^\infty$ that the reflections in its roots generate. 
The approximation results appear in the final section, though they can be read directly after \S\ref{HInfinity}.

Icosahedral symmetry and the Coxeter groups $H_2, H_3, H_4$ have continued to intrigue people from ancient times until the present day, where they are now familiar in such diverse places as carbon molecules, buckyballs, the capsids (outer shells) of viruses, Penrose tilings, and the aperiodic order of quasicrystals. There remains continuing interest in the mathematical world too, for instance, \cite{Koca} which investigates the subgroup structure of $H_4$ in a quaternionic context, and \cite{Dechant1} which explores some of the ways of making affine extensions of the $H_4$ root systems, following the success of affinization of the root systems of the simple Lie algebras. Affinization is accomplished by extending the Coxeter diagram, and might be thought of as based on using translations to extend root systems. In the present paper, although the root systems involved are defined intrinsically from the $H_3$ and $H_4$ root systems and the Galois conjugation $\tau\mapsto \tau'$, effectively the result is to extend by using contractions. Indeed this contractive aspect is something seems worthy of fuller study. 

We might note here that there have been considerable advances in the understanding of Coxeter groups along with their associated geometrical meanings by using the Clifford algebras, see \cite{Dechant2} and its associated references. Although we have not used these ideas here, and indeed we require only a very limited number of tools, the Clifford algebra approach may offer new insights into the setting we are introducing.

\section{$H_{4}$} \label{H4}

The set $\Delta$ of roots of the root system of type $H_{4}$ can be presented in the standard form
of \eqref{120-roots}
as the union of three types of vectors in $\R^{4}$. We let $K=K_{0}\cup K_{1}$ 
denote the set of roots of types {\bf R[0,0]} and {\bf R[1,0]},
which form a crystallographic root system of type $D_{4}$, and 
$\dot\Delta$ for the set of roots of type {\bf R[1,1]}. Of course these distinctions are not intrinsic to the root system itself, but
only to our coordinatization of it. However this distinction plays a fundamental role in what is to follow.

Along with $\Delta$ we have its conjugate $\Delta'$ which is the set of conjugates of $\Delta$, and the corresponding
sets $K$ and $\dot\Delta'$ (it being irrelevant whether or not the dot operators are applied before or after conjugation).
$K =\Delta \cap \Delta'$.

The reflections in the roots
$\Delta$ generate the group we call $H_4$. It is a Coxeter group of type $H_{4}$. Similarly we have $H_4'$ generated by
the reflections in the roots of $\Delta'$,
also of type $H_{4}$. We are primarily interested in the group $H^{\infty} $ which is generated by
$H_4$ and $H_4'$ together. The reflections given by $K$, the roots in common to both systems, generate a subgroup $D$ of type $D_{4}$ 
For example 
\begin{equation}\label{D4Base}
(1,0,0,0),\;\frac{1}{2}(-1,-1,-1,\pm 1),\; (0,1,0,0), \;(0,0,1,0),
\end{equation}
is a base (with either choice of the sign).

We let $\cS =\cS[4]$ denote the unit sphere in $4$-space. Both $\Delta$ and $\Delta'$ are in $\cS$, and so too
is the set $\Sigma = \Sigma[4]:=H^{\infty}(\Delta \cup \Delta')$. It is the smallest subset containing $\Delta \cup\Delta'$ and closed under
its own reflections (i.e. if $a \in \Sigma$ and $r_{a}$ is the reflection in $a$ then
$r_{a}(\Sigma) = \Sigma$). 

We note the important fact that 
\begin{equation}\label{cS1}
\frac{1}{2}\Ztau^{4} \cap \cS = \Delta \cup \Delta' \,.
\end{equation}

At the start we use only the algebraic consequences of reflections applied in the context of $\Z[\tau]^{4}$ and $\R^{4}$, but later we will interpret everything in the real quaternions, where the elements that we
are discussing can be interpreted as elements of \mbox{\rm{SU}$(2)$}.

It is useful to keep in mind the basic
facts about $\tau$:
\begin{equation}	
\tau + \tau' =1, \quad \tau\tau' = -1, \quad \tau^{2}= \tau +1, \quad \tau'^{2} = \tau' +1 \,.
\end{equation}

The quotient ring $\F_{4}:= \Ztau/2\Ztau$ is the Galois field of $4$ elements. Let $a\mapsto \overline a$ 
denote the natural homomorphism of $\Ztau \longrightarrow \F_{4}$. The elements of $\F_{4}$
are $\overline 0, \overline 1, \overline\tau,  \overline{\tau'}$. We extend the meaning of the map 
$\overline{\,\cdot \,}$
to its
component-wise version $\Ztau^{4} \longrightarrow \F_{4}^{4}$.

Inside the $\F_{4}$-vector space $\F_{4}^{4}$, define $A$ to be the subspace spanned by the elements of
$2\Delta$ taken modulo $2\Ztau^{4}$, i.e., $\overline{2\Delta}$. This is a $2$-dimensional space with basis
$(\overline 1, \overline 1,\overline 1,\overline 1)$ and $(\overline 0, \overline 1, \overline \tau', \overline \tau)$.  
We define $\mathring A:= A\backslash\{0\}$ and $\dot A := \overline{2\dot\Delta}$.
Note the cardinalities: $\card(A) =16$, $\card (\mathring  A)=15$, and $\card (\dot A)=12$.

Similarly we have
$A'$ using $2\Delta'$ and the corresponding subset $\dot A' = \overline{2\dot\Delta'}$. 
As we have already mentioned, and this sentence suggests, and will be true throughout, everything we say will come in two versions, which
are interchanged by conjugation. Henceforth we will usually only state and give proofs for one of the versions,
understanding that the other will be equally true. 

With the obvious dot product on $\F_{4}$ we find that $A$ is totally isotropic. This is just another way of saying that 
for all $\alpha, \beta \in K_{1}\cup\dot\Delta$, $2\alpha. 2\beta \in 2\Ztau$. In addition, if $\alpha$ and $\beta$ are linearly dependent, then $2\alpha. 2\beta \in 4\Ztau$. This simply follows
from $(0,1,\tau',\tau).(0,1,\tau',\tau)=4 = (1,1,1,1).(1,1,1,1)$.

We have parallel statements for $A'$ and $\dot A'$. However, for any
$a \in \dot A$ and $b \in \dot A'$
we always have $a.b \ne \overline 0$. 
\begin{tvr} \label{isotropy} \
\begin{itemize}
\item [{\rm (i)}] $A.A = \overline 0$. In particular, for all $x,y \in \Delta$,  $2x.y \in \Z[\tau]$;
\item [{\rm (ii)}] for all $a \in \dot A$ and $b \in \dot A'$, $a.b \ne \overline 0$. In particular if 
$x\in \dot\Delta$ and $y \in \dot \Delta'$ then $2x.y \notin \Z[\tau]$.
\end{itemize} \qed
\end{tvr}

\section{$H^{\infty}$ as a reflection group}\label{HInfinity}
\begin{figure}\centering
\includegraphics[scale=0.5]{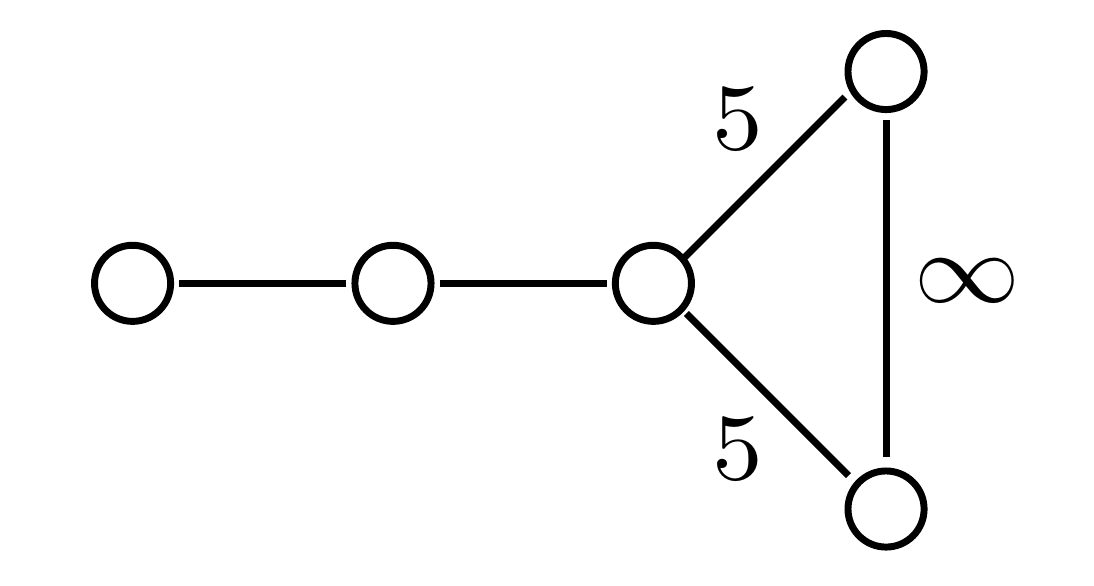}
\caption{This is the Coxeter diagram for the pairwise relationships between the reflection generators
of $H^{\infty}[4]$. The bonds between any two nodes indicate the orders of the product of the reflections
corresponding to them, with the convention that unmarked bonds indicate
order $3$ and bonds between commuting reflections are omitted completely. The proof establishing 
the bond marked with $\infty$ is found in Cor.~\ref{infiniteOrder}. We shall see that in fact $H^\infty[4]$ is a proper factor of the Coxeter group
of this diagram. The corresponding Coxeter diagram for the relationships between the generators of
$H^\infty[3]$ is obtained by deleting the first node and its attached bond. }
\label{HinftyDiagram}
\end{figure}

The set of vectors
\begin{equation} \label{Coxeter generators}
\{a_{1},a_{2},a_{3},a_{4},a_{5}\}= \{(1, 0, 0, 0), -\frac{1}{2}(1, 1, 1, 1), (0, 0, 1, 0),
\frac{1}{2} (0, -1, \tau', \tau), \frac{1}{2}(0, 1, -\tau, -\tau')\}\end{equation}
define reflections that correspond to the Coxeter diagram of Fig.~\ref{HinftyDiagram}
and these reflections generate $H^{\infty}= H^\infty[4]$. Indeed they contain a set of generators
for both of the Coxeter groups $H_{4}$ and $H_{4}'$. The associated Gram matrix 
with entries $a_{i}.a_{j}$ is positive semidefinite with null vector 
$(3- \tau, 6-2 \tau, 9-3\tau, -6+4\tau,2\tau)$. Certainly $H^{\infty}$ is a
factor of the corresponding Coxeter group, and later on we shall see that it is a proper factor,
see Prop.~\ref{HinfityRelations}.

\begin{tvr} \label{squaresProp}
If $x=(x_{1},x_{2},x_{3},x_{4}) \in \Ztau^{4}$ satisfies
\[x_{1}^{2 }+ x_{2}^{2 }+x_{3}^{2 }+x_{4}^{2 } \equiv 0 \mod 4\Ztau\]
then one of the three following cases prevails:
\begin{equation}
\begin{cases}
\mbox{all the $x_{i}$ are congruent to $0$ modulo $2\Ztau$,}\\
\mbox{ all the $x_{i}$ are congruent modulo $2\Ztau$
to each other, but not to $0$,}\\
\mbox{$(x_{1},x_{2},x_{3},x_{4}) $ is congruent modulo $2\Ztau^{4}$ to a permutation of the
components of  $(0,1,\tau',\tau)$. }
\end{cases}
\end{equation}
In the language and notation introduced above, these three cases correspond to 
\begin{equation} \label{threeCases}
\begin{cases}
\overline x = (\overline 0, \overline 0, \overline 0, \overline 0)\\
\overline x \in \mathring A \backslash (\dot A \cup \dot A')= \mathring A' \backslash (\dot A \cup \dot A')\\
\overline x \in \dot A  \mbox{  {\rm or} } \overline x \in \dot A'\,.
\end{cases}
\end{equation}

\end{tvr}
\noindent{Proof:} Modulo $2\Ztau$ each $x_{i} \equiv 0,1,\tau', \tau$. Their squares
are $\{0,1,1+\tau', 1+\tau\}$ modulo $4\Ztau$.
A quick check shows that the sum of four such elements can be zero modulo $4\Ztau$ if and only if one of the three conditions
of the thesis is true.  \qed

\medskip

Examples of all three types occur  when we double each of the roots, i.e. when we look at $2(\Delta \cup \Delta')$.

\begin{cor}\label{3squaresProp} The three-dimensional version of Prop.~\ref{squaresProp} reads:
if $x = (x_1,x_2,x_3) \in \Ztau^3$ satisfies
\[x_{1}^{2 }+ x_{2}^{2 }+x_{3}^{2 } \equiv 0 \mod 4\Ztau\] 
then one of the two following cases prevails:
\begin{equation}
\begin{cases}
\mbox{all the $x_{i}$ are congruent to $0$ modulo $2\Ztau$,}\\
\mbox{$(x_{1},x_{2},x_{3})$ is congruent modulo $2\Ztau^{3}$ to a permutation of the
components of  $(1,\tau',\tau)$. }
\end{cases}
\end{equation} \qed
\end{cor}

\medskip

For $n=1,2,3 \dots$ define the spherical sets
 \[ \cS_{n}:= \left \{x \in \frac{1}{2^{n}} \Ztau^{4}: x\in \cS,\;{\mbox {\rm and}}\;\; \overline {2^{n}x} \in \mathring A\right\}\, , \]
and 
\[\dot\cS_{n}:= \left \{x \in \frac{1}{2^{n}} \Ztau^{4}: x\in \cS,\;{\mbox {\rm and}}\; \;\overline{2^{n} x} \in \dot A\right\} \,;\]
and similarly the conjugated versions using $A'$. 
We call $n$ here the {\em level} of the spherical set. We define
\[K_n = K_n':=\cS_n \cap \cS_n' \,,\]
so $\cS_n= \dot\cS_n \cup K_n$, and similarly $\cS_n'= \dot\cS_n' \cup K_n$.

To these we add \[\cS_{0} = \cS_{0}':= K_{0}\,,\]
the set of roots of type {\bf R[0,0]}.

Here are two consequences of Prop.~\ref{squaresProp}:

\begin{tvr}
$\cS_{1} = K_{1}\cup\dot\Delta$, and 
$\dot \cS_{1}= \dot\Delta$.
\end{tvr}
\noindent{ Proof:} Let $y \in \cS_{1}$. Then $\overline{2y} \in \mathring A$ and $2y=(y_{1},y_{2},y_{3},y_{4}) \in \Ztau^{4}$ satisfies
$y_{1}^{2 }+ y_{2}^{2 }+y_{3}^{2 }+y_{4}^{2 } = 4$. Using Prop.~\ref{squaresProp}, one can hunt through all 
solutions to see that only elements of $\Delta \cup \Delta'$ can satisfy this, and then restrict to $\dot A$ filter them out
to $\dot\Delta$.  \qed

\begin{tvr} \label{unionSs}
For $n\ge 1$,
\[\frac{1}{2^{n}}\Ztau^{4} \cap \cS = \bigcup_{m\le n}(\cS_{m} \cup \cS_{m}') \,.\]
Furthermore, $x\in \frac{1}{2^{n}}\Ztau^{4}\cap \cS $ is in $\cS_{n} \cup \cS_{n}'$ if and only if 
$\overline{2^{n}x} \ne \{\overline 0,\overline 0,\overline 0,\overline 0\}$.
\end{tvr}  
\noindent{ Proof:} Let $x \in \frac{1}{2^{n}}\Ztau^{4} \cap \cS$, and let $m \le n$ be minimal so that 
$2^{m}x \in \Ztau^{4}$. If $m=0$ then $x\in K_{0}$. If $m\ge1$
 then $2^{m}x.2^{m}x = 4^{m}$ and by Prop.~\ref{squaresProp}, $\overline {2^{m}x} \in \mathring A\cup \mathring A'$, see \eqref{threeCases},
 from which $x\in \cS_{m} \cup \cS_{m}'$. The reverse inclusion comes from the definitions, and the last line is clear.
\qed

\begin{tvr}\label{Delta'OnSn}

For all $\alpha' \in \dot\Delta'$, $r_{\alpha'}\dot\cS_{n} \subset \dot\cS_{n+1}'$.
\end{tvr}
\noindent{Proof:} Let $x \in \dot\cS_{n}$, so by definition $2^{n}x \in \Ztau^{4}$ and $\overline{2^{n}x} \in \dot A$. We have
$ 2^{n} r_{\alpha'}x = 2^{n}x - (2^{n}x.2\alpha')\alpha'$. The elements $\overline{2^{n}x}$ and $\overline{2\alpha'}$
are in $\dot A$ and $\dot A'$ respectively, so in all cases $\overline{(2^{n}x.2\alpha')} \ne \overline 0$.
Since $\alpha'$ is not in $\Ztau^{4}$, but only in $(1/2)\Ztau^{4}$, 
$2^{n} r_{\alpha'}x= 2^{n}x - (2^{n}x.2\alpha')\alpha' \notin \Ztau^{4}$. However 
$2^{n+1} r_{\alpha'}x= 2^{n+1}x - (2^{n}x.2\alpha')2\alpha' \in \Ztau^{4}$ and furthermore $\overline{2^{n+1} r_{\alpha'}x}
=(\overline{2^{n}x.2\alpha'}) \overline{2\alpha'} \in \dot A'$\,.
Thus $r_{\alpha'}(x) \in \dot \cS_{n+1}'$, which is what we wished to prove. \qed

\medskip

\begin{cor} \label{infiniteOrder}
For all $\alpha \in \dot \Delta$ and for all $\beta' \in \dot\Delta'$, $r_{\alpha}r_{\beta'}$ has infinite order.
\end{cor}
\noindent{Proof:} Let $x\in \dot\cS_{1}$. Then $r_{\beta'}x \in \dot \cS_{2}'$, $ r_{\alpha}r_{\beta'}x \in \dot\cS_{3}$, 
$ r_{\beta'}r_{\alpha}r_{\beta'}x \in \dot\cS_{4}$
and so on. \qed
\bigskip

Let $H^\infty_+[4]$ denote the commutator subgroup of $H^\infty[4]$.
\begin{tvr}\label{denseness}\
\begin{itemize}
\item[{\rm(i)}] $H^\infty_+[4]$ is dense in $SO(4)$;
\item[{\rm(ii)}] $H^\infty[4]$ is dense in $O(4)$;
\item[{\rm(iii)}] $\Sigma[4]$ is dense in $\cS[4]$.
\end{itemize}
\end{tvr}
\noindent {Proof:} 
With $\alpha := a_{4}=(1/2)(0,-1,\tau', \tau) \in \dot\Delta$ and $\beta': = a_{5}=
(1/2)(0,1,-\tau, -\tau') \in \dot\Delta'$, $r_{\alpha} r_{\beta'}$ is a rotation in the plane
$P$ orthogonal to the two roots $(1,0,0,0)$ and $\frac{1}{2}(1,1,1,1)$, and by Cor.~\ref{infiniteOrder} it has infinite
order. Let $G_1\subset H^\infty_+[4]$ be the infinite cyclic group generated by $(r_{\alpha} r_{\beta'})^2= r_{\alpha} r_{\beta'}r_{\alpha} r_{\beta'} \in H^\infty_+[4]$. The orbit of $\alpha$ under $G_1$ generated is infinite and hence is
dense on the circle on the sphere $\cS \cap P$.  In particular $G_1$ has rotations through angles as small
as we please. 

Now also $\gamma':=r_{a_3}(a_5) = (1/2)(0,1,\tau, -\tau') \in \dot\Delta'$, so
$r_{\alpha} r_{\gamma'}$ is a rotation of infinite order in a different plane orthogonal to $(1,0,0,0)$ and we get
a second group $G_2 \subset H^\infty_+[4]$ with arbitrarily small rotations. These two groups generate a subgroup of $H^\infty_+[4]$ that
is a dense subgroup of the group of $SO(3)_1$ of rotations of $(0,\R,\R,\R)$. Since $H^\infty_+[4]$ can certainly 
map $(0,\R,\R,\R)$ to $(\R,\R,\R,0)$, we see that it also contains a subgroup that is a dense subgroup
of the $SO(3)_4$ of rotations of this space. Since $SO(3)_1$ and $SO(3)_4$ generate the entire rotation group
$SO(4)$, we see that $H^\infty_+[4]$ contains a dense subgroup of this space. 

Now (i),(ii),(iii) are all clear. \qed

\begin{tvr}\label{DeltaOnSn}
Let $n\ge 2$. Then
\begin{itemize}
\item[{\rm (i)}]  for all $\alpha \in \dot\Delta$,  $r_{\alpha}(S_{n})\subset \cS_{n}\cup \dot\cS_{n-1}'$;
\item[{\rm (ii)}] $r_{\alpha}(\cS_{n}\backslash \dot\cS_{n})\subset \cS_{n}$ for all $\alpha \in \dot\Delta$, 
and for each $x \in \cS_{n}\backslash \dot\cS_{n}$ there  exists $\alpha \in \dot\Delta$ for which
$r_\alpha(x) \in \dot\cS_n$;
\item[{\rm (iii)}] for each $x \in \dot\cS_{n}$ there are three reflections of the form $r_{\alpha}$ where $\alpha \in \dot\Delta$, for which 
$r_{\alpha}x \in\dot\cS_{n-1}'$.
\end{itemize}
\end{tvr}
\noindent{ Proof:} (i). We start with the slightly weaker assumption that $n\ge 1$. Let $x\in \cS_{n}$ and $\alpha \in \dot\Delta$.
\begin{equation}\label{reflFormula}
 2^{n} r_{\alpha}x = 2^{n}x - (2^{n}x.2\alpha)\alpha \,.
 \end{equation}

Let $\overline{2^{n}x} = b \in \mathring A$. Since $\overline{2^{n}x}$ and $\overline{2\alpha}$ are elements of $A$, $(2^{n}x.2\alpha) \in
2\Ztau^{4}$. Let $\overline{(2^{n}x.2\alpha)\alpha}=a \in A$. As long as $ a\ne b$, 
$\overline{ 2^{n} r_{\alpha}x } = b-a \in \mathring A$, so $r_{\alpha}x \in \cS_{n}$.
If $a=b$ then $2^{n} r_{\alpha}x \in 2\Ztau^{4}$, so $2^{n-1} r_{\alpha}x \in \Ztau^{4}$. 
By Prop.~\ref{unionSs}, $r_{\alpha}x$ is either in $\cS_{m}$ or $\cS_{m}'$ for some 
$m<n$.

Now whatever the case is, if $n\ge 1$ then we have shown that $r_{\alpha}(x)$ is either in $\cS_{n}$ or in some 
$\cS_m$ or $\cS_{m}'$ for $m<n$. Now assume that $n\ge 2$. Then if $m < n$ it is not possible that
 $r_{\alpha}(x) \in \cS_{m}$, for either
$m\ge1$ in which case by the same reasoning, $x= r_{\alpha}r_{\alpha}(x) $ is on some sphere of level at most $m$, contrary
to $x \in \cS_{n}$, or $m=0$ in which case $r_{\alpha}r_{\alpha}(x) $ has level $1$. In the case that $r_{\alpha}(x) \in \cS'_m$, then noting Prop.~\ref{Delta'OnSn} we see that actually $r_{\alpha}(x) \in \cS'_{n-1}$.
But in fact $r_{\alpha}(x) $ cannot be in $\cS_{n-1}'\backslash \dot\cS'_{n-1}$, since these elements do not change
level under the application of $r_{\alpha}$ and so $x= r_{\alpha}r_{\alpha}(x) $ would be impossible. Thus
$r_{\alpha}(x) \in \dot\cS'_{n-1}$
This proves part (i).

This last case, where the level drops under reflection, is particularly interesting and we wish to look more carefully how it can happen in (iii) below.

\medskip
(ii). Let $\alpha \in \dot\Delta$. In the notation of part (i) above, if $x\in \cS_n \backslash \dot\cS_n$ then $b:=\overline{2^n x}\in \F^\times_4 (\overline 1, \overline 1, \overline 1, \overline 1).$  Now if $a:= \overline{(2^{n}x.2\alpha)\alpha} \ne \overline 0$ then $a \in \dot A$, and then also $b-a \in \dot A$ and the resulting
$r_\alpha x \in \dot \cS_n$. However, if $a = \overline 0$ then $b -a = b$ and $r_\alpha(x)$ remains in
$\in \cS_n \backslash \dot\cS_n$. We need to show that we can choose $\alpha \in \dot\Delta$
for which the corresponding $a \ne \overline 0$. 
 
 For definiteness (the other cases work the same way) we assume that $b= (1,1,1,1)$ and write $2^n x =b + 2v$, where $v = (v_1,v_2,v_3,v_4) \in \Z[\tau]$. We begin
 with $\alpha =\frac{1}{2}(0,1,\tau', \tau)$, though ultimately we will cycle its last three components around. 
 Then the question is what does it mean for 
 \begin{equation}\label{badCongruence}
 a=\overline{(2^n x . 2\alpha) \alpha}= \overline 0
  \end{equation} 
 to hold, or quivalently, what does it mean for 
 $2^n x . 2\alpha= 2 + 2v.2\alpha \equiv 0 \mod 4\Z[\tau]$?

 The latter is equivalent to $1+ v.2\alpha \equiv 0
 \mod 2 \Z[\tau]$, or finally 
 \[ \overline{v_2} + \overline{\tau' v_3} + \overline{\tau} \overline{v_4} = \overline 1\,.\] 
 This may happen, but if we cycle the last three terms of $\alpha$ around then we get the same equation
 but with the coefficients $v_2, v_3,v_4$ likewise cycled around. Adding all three of the equations together
 and using the fact that $\overline 1 +\overline {\tau'} +\overline{\tau} = \overline 0$ we get the contradiction
 $\overline 0 = \overline 1$. Thus, for at least one of these three possibilities for $\alpha$, equation
 \eqref{badCongruence} fails, and this version of $\alpha$ produces $r_\alpha(x) \in \dot \cS_n$.
 
\medskip
 
(iii). Suppose that $x\in \dot\cS_{n}$ and let $\alpha \in \dot A$. What is required to have  $\overline{2^{n}r_{\alpha}x} =\overline 0$\,?
To see what happens here, it is easiest to fix one particular form that $\overline {2^{n}x}$ can take, and we 
choose $(\overline 0, \overline 1, \overline \tau', \overline{\tau})$. All other forms are even permutations of
this, and as will be seen, there is no loss in generality in assuming this form.  Now for any
$\alpha \in \dot\Delta$, $\overline{2\alpha}$ is some even permutation of 
$(\overline 0, \overline 1, \overline \tau', \overline{\tau})$. Referring to \eqref{reflFormula}, the question becomes, how can 
$\overline{(2^{n}x.2\alpha)\alpha} = a= (\overline 0, \overline 1, \overline \tau', \overline{\tau})$\,?
Looking at the $0$ in the first coordinate of $\overline {2^{n}x}$, it is clear that $\alpha$ must
be one of the elements, $(0,\pm 1,\pm\tau', \pm\tau)$,  $(0,\pm\tau', \pm\tau,\pm 1)$, or
$(0,\pm\tau,\pm1,\pm\tau')$. Since an overall sign change in $\alpha$ makes no difference to the 
reflection defined by $\alpha$ nor to the equation we are trying to solve, we can replace
the $\pm \tau'$ in each of these by simply $\tau'$. Doing this we find there is exactly one solution 
in each of the three sets. 

For instance, taking the second type we have 
\begin{equation}\label{quadEq}
(2^{n}x.2\alpha)\alpha =
\frac{1}{2}\left((0, 1, \tau', \tau).(0,\tau', u\tau,v 1)\right)(0,\tau', u\tau,v 1)=
\begin{cases}
0(0,\tau', \tau, 1) \stackrel{\mod 2}{\longrightarrow }(\overline 0, \overline 0,\overline 0,\overline 0)\\
1(0,\tau', -\tau, 1) \stackrel{\mod 2}{\longrightarrow }(\overline 0, \overline {\tau'},\overline \tau,\overline {1})\\
-\tau(0,\tau', \tau,-1) \stackrel{\mod 2}{\longrightarrow} (\overline 0, \overline {1},\overline {\tau'},\overline \tau)\\
\tau'(0,\tau', -\tau,- 1) \stackrel{\mod 2}{\longrightarrow} (\overline 0, \overline \tau,\overline {1},\overline {\tau'})
\end{cases}
\end{equation}
as $u,v$ go through the various possibilities of sign, so $u=1$, $v=-1$ gives the value we are
looking for. In fact the other two solutions are just obtained by cyclically permuting the last three coordinates,
so the solutions in this case are 
\begin{equation} \label{3solutions}
(0,\tau',\tau,-1), \quad (0,\tau,-1,\tau'), \quad (0,-1,\tau',\tau) \,.
\end{equation}

The three solutions, and only these three, produce the desired effect that $2^{n}r_{\alpha}(x) \in 2\Ztau^{4}$,
and as a consequence, as we saw above,  $r_{\alpha}(x) \in \dot \cS_{n-1}'$. This concludes the proof of part (iii).
\qed

\medskip
It is worthwhile noting that the three vectors we have found in \eqref{3solutions} are all in the plane orthogonal to the 
roots $(1,0,0,0)$ and $\frac{1}{2}(1,1,1,1)$ and they lie at angles of $2\pi/3$ to each other. 
\medskip

Let 
\begin{equation} \label{DefOfDGN}
\cR:= \bigcup_{n=0}^{\infty}\frac{1}{2^{n}}\Ztau
      = \left\{\frac{a+b\tau}{2^n}\,:\, a,b \in \Z, n\in \Z_{\ge 0}\right\} = \Z[\frac{1}{2},\tau]\,.
\end{equation}

This is a ring which we shall call the ring of {\em dyadic golden numbers}.
\begin{tvr}\label{fullness}
\[\Sigma[4] = \bigcup_{n=0}^{\infty} (\cS_{n}\cup \cS_{n}') = \cS[4] \cap \cR^{4} \,.\]
\end{tvr}
\noindent
Proof: By definition \eqref{defSigma}, $\Sigma= H^{\infty}(\Delta \cup \Delta'$). Let $\Sigma_{n} := \cS_{n} \cup \cS_{n}'$ for all $n=0,1,2,\dots$. 
We have $\Sigma_{0} \cup \Sigma_{1} = \Delta \cup \Delta' \subset \Sigma$.
Now proceed by induction and assume that we have shown that $\Sigma_{n} \subset \Sigma$. Then from Prop.~\ref{Delta'OnSn}
and Prop.~\ref{DeltaOnSn}(iii),
$\dot\cS_{n+1}\cup \,\dot\cS_{n+1}' \subset \Sigma$. By Prop.~\ref{DeltaOnSn}(ii) we see
that all of $\Sigma_{n+1} \subset \Sigma$. Thus 
\[\bigcup_{n=0}^{\infty} (\cS_{n}\cup \cS_{n}') \subset \Sigma\,.\]
However Prop.~\ref{Delta'OnSn} and Prop.~\ref{DeltaOnSn} show that the left-hand side here is closed
under the action of the reflections in the roots $\dot\Delta \cup \dot\Delta'$, and these
reflections generate all of $H^{\infty}$.
This gives the reverse inclusion.  For the final equality see Prop.~\ref{unionSs}.\qed

\begin{tvr}
For all $n\ge1$, $\card(\dot\cS_{n}) = \card(\dot\cS_{n}') = 3.(2)^{4n+1}$.
\end{tvr}
\noindent
\noindent Proof: $\dot\cS_{1}$ consists of the $96$ roots of $\dot \Delta$, and similarly for $\dot\cS_{1}'$. 
Now proceed by induction. We simply have to note that for each element of $\dot\cS_{n}$ there
are $48$ reflections arising from $\dot\Delta'$ that will map it into $48$ images in $\dot\cS_{n+1}$, 
see Prop.~\ref{Delta'OnSn}.
However, according to Prop.~\ref{DeltaOnSn}(c) each of these elements will be produced exactly three times
as we use all the reflections available from $\dot\Delta'$. Thus $\card \dot\cS_{n+1}'= 16\,\card\dot\cS_{n}$,
and likewise for $\card \dot\cS_{n+1}$. This completes the induction step. \qed

\medskip
The consequence of Prop.~\ref{DeltaOnSn} and Prop.~\ref{Delta'OnSn} is that to go deeper into the structure
of $\Sigma$ we need to apply reflections alternately from $\dot\Delta$ and $\dot\Delta'$.
It is this path that leads to finer and finer discretization of the unitary group.

\medskip

\medskip

\begin{figure}
\centering
\includegraphics[scale= 0.3]{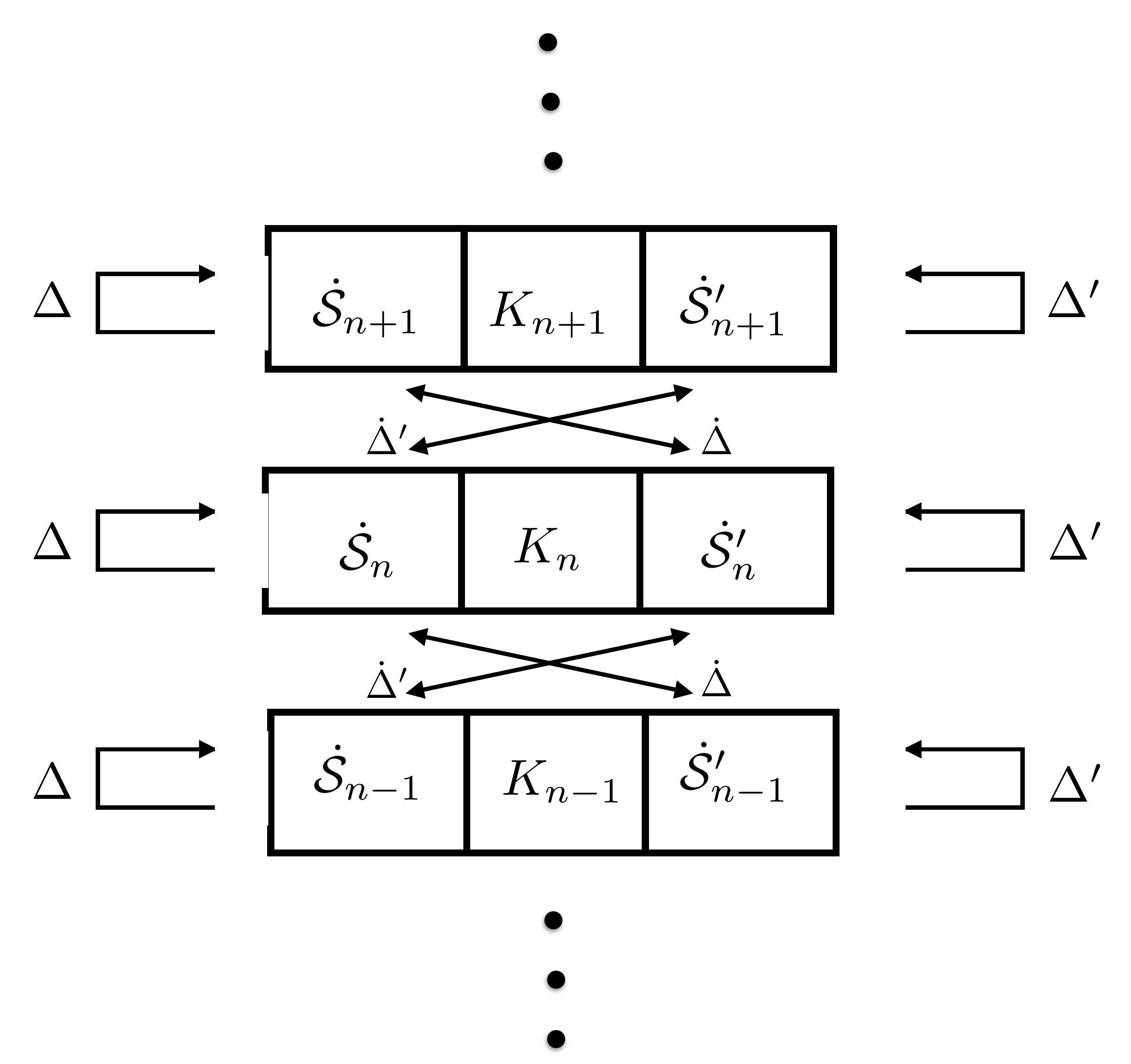}
\caption{The action of reflections from different types of roots is indicated diagramatically.
The intersection  $\cS_{n} \cap \cS_{n}'=:K_n$ is the common `crystallographic' part of 
the two spherical sets.  The lateral actions of $\Delta$ apply to $\cS_{n}$ but not
to $\cS_{n}'$, and similarly with $\Delta'$.}
\label{SnFlow}
\end{figure}

The situation that we have just described is pictorially represented in Fig.~\ref{SnFlow}.

\section{$H_{3}$}\label{H3}

The four-dimensional theory that we have outlined so far has a parallel version in three dimensions around the 
Coxeter group $H_3$. Most of the features that we have described in the four-dimensional case appear
again here, but now notationally distinguished with the additional suffix `$[3]$'. When we are dealing with both
settings together, we shall also use the suffix `$[4]$' for the four-dimensional situation in order to avoid confusion. 

The unit sphere in $(1/2)\Ztau^{3}$ consists of $54$ points, namely the $6$ vectors $(\pm 1,0,0)$ and all permutations of its coordinates, and the $48$ points
$(1/2)(\pm 1, \pm\tau', \pm \tau)$ and all permutations of these coordinates. If in the second set of points we restrict them by allowing
only even (= cyclic) permutations then they provide $24$ points which together with the original $6$, produce a set of $30$ points which may be viewed as the 
non-crystallographic root system $\Delta[3]$ of type $H_{3}$. They are the vertices of the
icosadodecahedron. The reflection group generated by the 
$15$ distinct reflections in these points, now called roots, that is the reflections in the mirrors (planes) through the origin
orthogonal to these roots, is the icosahedral group $H_3$ of order $120$. 

As a base for the root system of type $H_3$ we can use
the vectors $(1/2)(1,\tau',\tau)$, $(1/2)(1,-\tau',-\tau)$, $(1/2)(\tau',\tau,1)$.

Instead of the even permutations of $(1/2)(\pm 1, \pm\tau', \pm \tau)$ we can conjugate them, thus taking only 
all the even permutations of $(1/2)(\pm 1, \pm\tau, \pm \tau')$ (or equivalently we take all the odd permutations)
and do the same thing. We then end up with $\Delta[3]'$, the set of conjugates of the elements of $\Delta[3]$.
Again we have a root system of type $H_{3}$ and again the reflections in these roots generates a copy
of the icosahedral group. 

\smallskip
We write $\dot\Delta[3]$ for the non-rational roots of $\Delta[3]$, and similarly
for $\dot\Delta[3]'$.

\smallskip
Let $H^{\infty}[3]$ denote the group generated by the reflections in 
$\Delta[3] \cup \Delta[3]'$ and let $\Sigma[3]$
be the orbit of $\Delta[3] \cup \Delta[3]'$ under $H^{\infty}[3]$.  We shall see that $\Sigma[3]$ is a dense set of points on the unit $2$-sphere $\cS[3]$.

We put this into the $4$-dimensional setting in which we have discussed $H_{4}$ by considering the $3$-dimensional space
here as the subspace $V:=(0,\R,\R,\R)$. Then it is 
 clear that  $\Delta[3] = \Delta[4]\cap V$ and $\Delta[3]'= \Delta[4]'\cap V$, and indeed all of $\Sigma[3]$, lies in $ V\cap \Sigma[4]$. In fact we shall see that 
$\Sigma[3] = V \cap \Sigma[4]$.

Later on we shall put everything into the setting of the algebra of quaternions $\HH$, at which point 
$V$ will be come the subspace of pure quaternions.

  As before, let $\F_{4}:= \Ztau/2\Ztau$ and let $x\mapsto \overline x$ 
denote the natural homomorphism of $\Ztau \longrightarrow \F_{4}$. In place of $\F_{4}^4$ 
we will now wish to use $\F_{4}^3$. It is often convenient to identify this with the  subspace $(0,\F_{4},\F_{4},\F_{4})$ 
of $\F_{4}^4$, and we will do this when it is useful to do so. 

 Using this we define the analogues of $A[4]$ and $A[4]'$, namely
 \begin{align}
A[3]&:= A[4] \cap V \,,\nonumber\\
 A[3]'&:= A[4]' \cap V \,
\end{align}
which will play the same role as $A$ and $A'$ in the four-dimensional case. 
Written as elements of $\F_{4}^3$ 
\[A[3]= \{(\overline 0, \overline 0,\overline 0),  (1,\overline{\tau'}, \overline{\tau}), 
 (\overline{\tau'}, \overline{\tau},1), (\overline{\tau},1,\overline{\tau'})\}\,,\]
 and similarly for $A[3]'$.
Define
$\dot A[3]:= \dot A[4] \cap V = \mathring A[4] \cap V$, and similarly $\dot A[3]'$.

The dot product on $\F_{4}^{3}$ is defined by $\overline u. \overline v = \overline{u.v}$
 for all $u,v \in \Ztau^{3}$, which matches the dot product it has as a subspace of $\F_{4}^4$. From
 Prop.~\ref{isotropy} we have:
 \begin{tvr} \label{isotropy3} \
\begin{itemize}
\item [{\rm (i)}] $A[3].A[3] = \overline 0$. In particular, for all $\alpha ,\beta \in \Delta[3]$,  $2\alpha.\beta \in \Z[\tau]$;
\item [{\rm (ii)}] for all $a \in \dot A[3]$ and $b \in \dot A'[3]$, $a.b \ne \overline 0$. In particular if 
$\alpha\in \dot\Delta[3]$ and $\beta \in  \dot\Delta[3]'$ then $2\alpha.\beta \notin \Z[\tau]$.
\end{itemize} \qed
\end{tvr}
  
  \medskip

We now obtain very similar results to those of \S\ref{HInfinity}.  We mainly just list these, the proofs being basically the same. However, there is one significant difference that appears immediately---there are no longer two spherical sets
to deal with. For $n\ge 1$, we define

 \[ \cS_{n}[3]:= \left \{x \in \frac{1}{2^{n}} \Ztau^{3}: x\in \cS[3],\;{\mbox {\rm and}}\;\; \overline {2^{n}x} \in 
 \dot A[3]\right\}\, , \]
and similarly $\cS_{n}[3]'$ using $A'$. 
We call $n$ here the {\em level} of the spherical set. The sets $K_n$ and $K_n'$ have no relevant counterparts if 
$n>0$ since
\[\cS_n[3] \cap \cS_n[3]'  = \emptyset\,.\] 

To these we add \[K_0[3] :=\cS_{0}[3] = \cS_{0}[3]'=\{\pm 1,0,0), (0,\pm1,0),(0,0,\pm1)\}\,.\]

It is clear that $\cS[3] = \cS[4] \cap V$ and $\cS_n[3]= \cS_n[4] \cap V$, and parallel statements
for $\cS'$.

As a consequence of Cor.~\ref{3squaresProp}:
\begin{tvr} \label{unionSs3} \
\begin{itemize}
\item[{\rm (i)}] $\cS_{1}[3] = \dot\Delta[3]$;\\
\item[{\rm (ii)}] for all $n\ge 0$,
\[\frac{1}{2^{n}}\Ztau^{3} \cap \cS[3] = \bigcup_{m\le n}(\cS_{m}[3] \cup \cS_{m}'[3]) \,.\]
\end{itemize}
Furthermore, $x\in \frac{1}{2^{n}}\Ztau^{3}\cap \cS[3] $ is in $\cS_{n}[3] \cup \cS_{n}'[3]$ if and only if 
$\overline{2^{n}x} \ne \{\overline 0,\overline 0,\overline 0\}$.
\qed
\end{tvr}

 \medskip
 For each $\alpha \in \Delta[3]$, we also have $\alpha \in \Delta[4]$ and the reflection
 in $\alpha$ taken in $V$ is just the restriction of the reflection $r_\alpha$ taken in $\R^4$ restricted
 to $V$. We use the same symbol $r_\alpha$ for both, with context distinguishing them if necessary. 
 Corresponding primed notation will be used for $\Delta[3]'$. $H_{3}, H_{3}'$ denote
 the groups generated by the reflections of $\Delta[3]$ and $\Delta[3]'$ respectively. We may think of these
 groups as subgroups of $H_{4}$ and $H_{4}'$.
 
 From Prop.~\ref{Delta'OnSn} we have
\begin{tvr}\label{Delta'OnSn3}
For all $\alpha' \in \dot\Delta[3]'$, $r_{\alpha'}\cS_{n}[3] \subset \cS_{n+1}'[3]$. \qed
\end{tvr}  

\medskip

\begin{cor} \label{infiniteOrder3}
For all $\alpha \in  \dot\Delta[3]$ and for all $\beta' \in \dot\Delta[3]'$, $r_{\alpha}r_{\beta'}$ has infinite order.
\qed
\end{cor}
 
 Following the arguments of Prop.~\ref{denseness} we have 

\begin{tvr}\label{denseness3}\
\begin{itemize}
\item[{\rm(i)}] $H^\infty_+[3]$ is dense in $SO(3)$;
\item[{\rm(ii)}] $H^\infty[3]$ is dense in $O(3)$;
\item[{\rm(iii)}] $\Sigma[3]$ is dense in $\cS[3]$.
\end{itemize}
 \qed
\end{tvr}

Define $V_{\cR} := V \cap \cR^4$. It is now rather easy to use the four-dimensional results of \S\ref{HInfinity} to 
deduce corresponding results about the three-dimensional situation. We state these without further comment.

\begin{tvr}\label{Sigma3FromSigma4}
\[\Sigma[3] = \bigcup_{n=0}^{\infty} (\cS_{n}[3]\cup \cS_{n}'[3]) = \cS[3] \cap V_{\cR} = \cS[4] \cap V_{\cR}
=\Sigma[4]\cap V_{\cR}= \Sigma[4]\cap V\,.\]
\end{tvr}

\begin{tvr}
For all $n\ge1$, $\card(\cS_{n}[3]) = \card(\cS_{n}'[3]) = 3.(2)^{2n+1}$.
\end{tvr}
\noindent
\noindent Proof: $\cS_{1}[3]$ consists of the $24$ roots of $\Delta[3]$, and similarly for $\cS_{1}'[3]$. 
Now proceed by induction. We simply have to note that for each element of $\cS_{n}[3]$ there
are $12$ reflections arising from $\Delta[3]'$ that will map it into $12$ images in $\cS_{n+1}[3]$, 
see Prop.~\ref{Delta'OnSn}.
However, according to Prop.~\ref{DeltaOnSn}(iii) (which is actually set up so as to make the corresponding three dimensional
case rather apparent) each of these elements will be produced exactly three times
as we use all the reflections available from $\Delta[3]'$. Thus $\card \cS_{n+1}'[3]= 4\,\card\cS_{n}[3]$,
and likewise for $\card \cS_{n+1}[3]$. This completes the induction step. \qed

\section{How $\Sigma[4]$ acts on $\Sigma[3]$}

\begin{tvr}\label{transitive} \

\begin{itemize}
\item[{\rm (i)}] $H^{\infty}[4]$ acts transitively on $\Sigma[4]$\,;
\item[{\rm (ii)}] $H^{\infty}[3]$ acts transitively on $\Sigma[3]$\,;
\item[{\rm (iii)}] For each $\alpha \in\Sigma[4]$
the set of roots {\rm (}elements of $\Sigma[4]${\rm )} orthogonal to $\alpha$ is a translate of $\Sigma[3]$
by an element of $H^{\infty}[4]$\,; 
\item[{\rm (iv)}] For each $\alpha \in\Sigma[4]$ the stabilizer of $\alpha$ in $H^{\infty}[4]$ contains a subgroup 
conjugate to $H^{\infty}[3]$.
\end{itemize}
\end{tvr}
\noindent
{Proof:} (i) It suffices to show that $H^{\infty}[4] (1,0,0,0)$ contains the basic roots of \eqref{Coxeter generators}.
This happens because in the Coxeter diagram of Fig.~\ref{HinftyDiagram}, $(1,0,0,0)$ is connected to the four others  by chains of bonds that are all odd (namely labelled with $3$ or $5$). The argument is easily explained by doing it for the case of a bond labelled $3$ between two roots $\alpha, \beta$.
Then $r_{\alpha}\,r_{\beta}\,r_{\alpha}\,r_{\beta}\,r_{\alpha} = r_{\beta}$. But the left side is also
$r_{r_{\alpha}\,r_{\beta}(\alpha)}$, showing that $r_{\alpha}\,r_{\beta}(\alpha) = \pm \beta$. If the sign is wrong
we can use one more reflection $r_{\beta}$. 
The proof of (ii) is similar.

(iii) The stabilizing reflections $r_{a}$ of $(1,0,0,0)$, where $a\in \Sigma[4]$,  consists of all elements of $\Sigma[4] \cap V= \Sigma[3]$, see Prop.~\ref{Sigma3FromSigma4}. Parts (iii) and  (iv) follow directly from this. \qed 

\medskip

Instead of thinking of $\Delta[4]$ and $\Delta[4]'$ simply as subsets of $\R^{4}$, we can view
them as subsets of the quaternion algebra $\mathbb H = \R1 + \R i +\R j +\R k$.  We will identify
the subspace  $V$ of $\R^4$ that we have been using with the subspace
of pure quaternions: $V= \R i +\R j +\R k$. Then $\Delta[3] = \Delta[4]\cap V \subset \HH$
and similarly for $\Delta[3]'$. Generally we use this interpretation henceforth. 

We equip $\mathbb H$ with the
usual conjugation $u \mapsto \widetilde u$ which changes the signs of the last three components of each vector.\footnote{Usually the conjugation is
designated by the $\overline{\phantom{w}}$, but that would cause confusion with our earlier notation.} 
The usual dot product of $\R^{4}$
is the standard one for $\mathbb H$, and can be expressed in terms of the quaternionic multiplication by
\[x.y = \frac{1}{2} (x\widetilde y + y \widetilde x) \,.\]

As is well known, the unit sphere 
$\cS$ of $\mathbb H$, that is the set of vectors $a\in \mathbb H$ satisfying 
$\widetilde a a      = 1$, can be identified with $SU(2)$. Explicitly 
$\cS$ can be interpreted as the group $SU(2)$ by the usual matrix representation of $\mathbb H$ using 
\begin{equation} \label{matrixRepOfQuaternions}
1 \Leftrightarrow \begin{bmatrix} 1 & 0\\0&1 \end{bmatrix},\quad
i \Leftrightarrow \begin{bmatrix} \sqrt{-1} & 0\\0&-\sqrt{-1} \end{bmatrix} ,\quad
j \Leftrightarrow \begin{bmatrix} 0 & 1\\-1&0 \end{bmatrix},\quad
k \Leftrightarrow \begin{bmatrix} 0 & \sqrt{-1}\\ \sqrt{-1}&0 \end{bmatrix} \,.
\end{equation}

Returning to the ring \eqref{DefOfDGN}, we define $SU(2, \cR)$ to be the set
of all unitary matrices with coefficients in $\cR + i \cR$. Using \eqref{matrixRepOfQuaternions}, $SU(2,\cR)$
identifies with the set of all quaternions with coefficients in $\cR$ with norm $1$, i.e. with $\cR^{4}\cap \cS[4]$.

$SU(2,\cR)$ is a subgroup of $SU(2)$, and it is $\cS[4] \cap \cR^{4} = \Sigma[4]$ by Prop.~\ref{fullness}.
This tells us that $\Sigma[4]$ has the natural structure of a group, namely $SU(2,\cR)$. Furthermore, since
$\Sigma[4]$ is dense in $\cS[4]$, $SU(2,\cR)$ is dense in $SU(2)$.

We recall that for any $a\in \cS$ and any $x\in \mathbb H$ we have
\begin{equation}\label{quaternionReflection}
r_a(x) = x - (2x.a) a = x - (x\tilde a + a\tilde x)a = - a \tilde x a \,.
\end{equation}
In particular this shows that reflection can be effected by multiplication. Since $\Sigma[4]$ is by 
definition the smallest set containing $\Delta \cup \Delta'$ and closed under its own reflections, it
follows that $\Sigma[4]$ is generated as a group by $\Delta \cup \Delta'$.

\medskip
Now we are viewing $V$ as ordinary three-dimensional space, and since $V = \R i +\R j +\R k$, so $x\in V$ if and only if
$\widetilde x = -x$.  
Then
for all $a\in \cS$ and all $x\in V$, 
\[ \widetilde{ axa^{-1} } = \widetilde{ax\widetilde{a}} = a (-x) \widetilde a = - axa^{-1}  \,, \]
showing that $axa^{-1} \in V$. This is the familiar way in which elements of $SU(2)$
turn into elements of $SO(3)$, so that $a\in SU(2)$ maps to the conjugation
$\gamma_a: x \mapsto axa^{-1}$ on $V$. This is a rotation, and in this way we obtain the well-known double cover mapping
\begin{equation}\label{SU2toSO3}
 SU(2) \longrightarrow SO(3) \,.
 \end{equation}
However, there is more. 
If $a \in \Sigma[4]$ then for all $x\in \Sigma[3]$, then $ax\tilde a \in \Sigma[3]$, since
if $a \in \frac{1}{2^{n}}\Ztau^{4}$ and
$x \in \frac{1}{2^{m}}\Ztau^{3}$ then $ax\tilde a \in \frac{1}{2^{m+2n}}\Ztau^{3} \cap \cS[3]\subset \Sigma[3]$.
Thus the group of rotations on the sphere in $3$-space induced by $\Sigma[4]$ actually acts as a 
group of symmetries of $\Sigma[3]$. In particular they all lie in $SO(3,\mathcal R)$. 

\begin{tvr} \label{SigmaAsAGroup}\
\begin{itemize}
\item[{\rm(i)}] $\Sigma[4]= \cR^{4} \cap \cS[4]$ is a group under quaternion multiplication
and is generated by $\Delta \cup \Delta'$;
\item[{\rm(ii)}]
$\Sigma[4]= SU(2,\cR)$ and  $SU(2,\cR)$ is dense in $SU(2)$;
\item[{\rm(iii)}] $\Sigma[4]$ acts naturally as a group of isometries on $\Sigma[3]$
via the mapping $\gamma_a:x \mapsto ax\tilde a$ for all $a\in \Sigma[4]= SU(2,\cR)$; 
\item[{\rm(iv)}] Under the mapping $\gamma$ in {\rm (iii)}, $SU(2,\cR)$ maps onto a subgroup of
index $2$ in $SO(3,\cR)$.
\end{itemize}
\qed
\end{tvr}

\noindent
Proof:  There remains only to prove (iv). Let 
$\rho:= (\frac{1}{\sqrt2},\frac{1}{\sqrt2} ,0,0)$. We begin by making the observation that 
\begin{align}\label{someSimpleRotations}
\gamma_i((x_1,x_2,x_3,x_4)) &= (x_1, x_2, -x_3, -x_4) \nonumber \\
\gamma_j (x_1,x_2,x_3,x_4))  &= (x_1, -x_2, x_3, -x_4)\\
\gamma_k((x_1,x_2,x_3,x_4))  &= (x_1, -x_2, -x_3, x_4) \nonumber \\
\gamma_\rho ((x_1,x_2,x_3,x_4)) &= (x_1, x_2, -x_4, x_3) \,. \nonumber
\end{align}
Since the mapping $SU(2)\longrightarrow SO(3)$ is $2:1$ and its kernel is $\pm 1 \in SU(2)$,
it is clear that the last of these rotations, $\gamma_\rho$, while is is in $SO(3,\cR)$, is not something that can arise from $SU(2,\cR)$. 
We shall see that this is the only obstruction that we have to deal with. 

Since $\Sigma[3] = \cS \cap V_\cR$,
$SO(3,\cR)$ maps $\Sigma[3]$ into itself.  Let $G:=\gamma_{SU(2,\cR)}$ denote the image
of $SU(2,\cR)$ in $SO(3,\cR)$ under the mapping $SU(2) \longrightarrow SO(3,\cR)$.
Note that if $a \in \Sigma[3]$ and
$x \in V$ then $\tilde x = -x$ and \eqref{quaternionReflection} reads 
\begin{equation}\label{reflectEqualsMinusRot}
r_a(x) = -a \tilde x a = axa = - a xa^{-1}\,.
\end{equation}

 Now suppose
that $R \in SO(3,\cR)$. We wish to make $R$ as simple as possible by multiplying $R$ on the left
by elements in $G$. Consider $R(0,1,0,0) \in \cS \cap V_\cR = \Sigma[3]$. 
Since $H^\infty[3]$ is transitive on $\Sigma[3]$, there exists $w \in H^\infty[3]$ with 
$wR=(0, 1,0,0)$. Writing $w= r_{a_k} \dots r_{a_1}$, \eqref{reflectEqualsMinusRot} shows that 
$(-1)^k w \in G$. This being the case we can assume at the outset
that $R(0,1,0,0) = (0, \pm 1, 0,0)$. However we see from \eqref{someSimpleRotations} that we can
alter this sign, so we can actually assume that $R(0,1,0,0) = (0, 1, 0,0)$. Now, the only elements of $\Sigma[3]$ that
have the form $(0,0,*,*)$ are $(0, 0,\pm1,0)$ and $(0,0,0,\pm1)$, and since $R$ is an orthogonal transformation
fixing $(1,0,0,0)$ and $(0,1,0,0)$, it can only permute these two pairs roots around. Of the four possible
rotations that are possible here, two of them require interchanging the two coordinate positions and changing one sign, and that requires $\gamma_\rho$. Thus, by composing $R$ with elements of $G$ we can 
reduce to one of the pair $\{1, \gamma_\rho\}$; and indeed, $\gamma_\rho \notin G$. This proves that $[SU(2):G]=2$. \qed 

\medskip
We note that $SO(2,\cR)$ is not dense in $SO(2)$, since the only $\cR$-points on the
unit circle are $(\pm1, 0)$ and $(0, \pm1)$. Likewise $U(2,\cR)$ (the subgroup of all
the elements of $U(2)$ whose coefficients are in $\cR + i \cR$) is not dense in $U(2)$---and for basically the 
very same reason: $\det$ on $U(2)$ is surjective onto the unit circle.

\begin{tvr}
The stabilizer of $(1,0,0,0)$ in $\langle H^\infty[4],\iota\rangle$ is 
$\langle H^{\infty}[3], \iota\rangle$ where $\iota$ is the permutation operator on the last two coordinates. 
\end{tvr}

 Here $H^\infty[3]$ is in its obvious form of working on the last three coordinates, that is on the space $V = \mathbb R i + \mathbb R j + \mathbb R k$ in the quaternions, and $\langle H^{\infty}[3], \iota\rangle$ is the group 
 generated by $H^{\infty}[3]$ and $\iota$. 
 \medskip
 
\noindent Proof: The first observation is that in $H_3$ the rotation of the last three coordinates exists (as a rotation of order $3$)
and also the changes of signs of the individual coordinates exist (as the reflections in these coordinates). The same goes for $H_3'$, of course. 
Now let $w \in H^\infty[4]$ be in the stabilizer of 
$z:= ( 1,0,0,0)$. Let $y:= (0,1,0,0)$ and let $x:= w(y)$.
Then $x\in V$ since $w$ is an isometry. 
Since $H^\infty[3]$ is transitive on $\Sigma[3]$, there is a $v\in H^\infty[3]$ with $v(x) = y = (0,1,0,0)$.
Now the element $vw$ fixes both $x$ and $y$. We can ask what $vw$ does to $\Sigma[4]$.
It must map the remaining two elements $(0,0,1,0)$ and $(0,0,0,1)$ into roots of the form $(0,0,*,*)$.
However the general form of roots shows that the only possibilities are  $(0,0,\pm1,0)$ and $(0,0,0,\pm1)$. The signs can be adjusted as we please by using an element from the Weyl group of $K$, but
$vw$ may still interchange the two points, in other words may be the involution $\iota$. We shall show in 
Prop.~\ref{orthoGroupRevealed} that $H^{\infty}[4]$ does not contain $\iota$. \qed

\section{$\Delta$ and $\Delta'$ as groups}

The two basic root systems we began with are $\Delta$ and $\Delta'$. The rational parts of these sets
is their intersection $K= \Delta \cap \Delta'$.

\begin{tvr}\label{DeltaGroup}
$K$, $\Delta$, and $\Delta'$ are groups {\rm(}under quaternion multiplication{\rm)} of
orders $24,120,120$ respectively.
\end{tvr}

This is a well-known fact, see for instance \cite{DV} \S20. But for the convenience of the reader we prove it 
here. 

\noindent
Proof:  For the rational root system $K$ it is clear that the product of any two elements of $K$
is a unit vector in the rational space $(1/2) \Z^{4}$, and so is still in $K$. Inversion is just quaternion
conjugation, and it preserves $K$. Thus $K$ is a group with order $24$.  

We shall prove that $\Delta$ is a group; the proof for $\Delta'$ then follows immediately. The roots of $\Delta$ are given in \eqref{120-roots}.  The first set of these consists, in quaternion notation, of the elements
$\pm 1, \pm i, \pm j,\pm k$. Multiplication of a quaternion $(x_{1}, x_{2},x_{3},x_{4})$ on the left by $i$ turns it into
$(-x_{2}, x_{1}, -x_{4}, x_{3})$. There are sign changes, but the important thing to observe is that it performs an even permutation of the coordinates, interchanging the first two coordinates and the last two coordinates. Likewise multiplication by $j$ or $k$ 
performs an even permutation of the coordinates (and makes two sign changes). Similar things happen if the multiplication 
by $i,j,k$ is on the right. In particular, left and right multiplications by $\pm 1, \pm i, \pm j,\pm k$ map 
$\Delta$ into itself.

Now apart from the elements $\pm 1, \pm i, \pm j,\pm k$ in $\Delta$, every other $x \in \Delta$ has the property
that $2x \in \Z[\tau]^4$ and $\overline{2x} \ne (\overline 0, \overline 0,\overline 0,\overline 0)$ in $\mathbb F_{4}^{4}$. In 
fact $\overline{2x}$ is either $(\overline 1, \overline 1,\overline 1,\overline 1)$ or an even permutation of 
$(\overline 0, \overline 1,\overline {\tau'},\overline \tau)$. This property actually characterizes elements of norm $1$
that are in $\Delta \backslash \{\pm 1, \pm i, \pm j,\pm k\}$. Similar remarks hold for $\Delta'$ where the roles of $\tau$ and $\tau'$ are interchanged.

Now let $x,y \in \Delta$. We want to prove that $xy \in \Delta$. By what we said above, we can assume neither is in the set
$\{\pm 1, \pm i, \pm j,\pm k \}$. We begin by showing that $2xy \in \Z[\tau]^4$. By assumption $2x$ and $2y$ are in $\Z[\tau]^4$. We know that even permutations and arbitrary sign changes on $y$ will leave it in $\Delta$. 
In particular $\tilde y \in\Delta$, and so we know
that $2x.\tilde y \in \Z[\tau]$, by Prop.~\ref{isotropy}. Now $x.\tilde y$ is the first coordinate in the quaternion $xy$, so we know
that that the first coordinate in $2xy$ is $2x.y \in \Z[\tau]$. 

If we look at the second coefficient of $xy$
we can see that it is $(x_{1},x_{2},x_{3},x_{4}).(y_{2},y_{1}, y_{4},-y_{3})$, so this shows
that the second coefficient of $2xy$ is in $\Z[\tau]$. Continuing this idea with the
third and fourth coefficients we end up with $2xy \in \Z[\tau]$. 

At this point, since $|xy|= |x|\,|y| = 1$ we know that $xy$ is an element of  $\Delta$ or $\Delta'$. However, it is not possible for $xy \in \dot\Delta'$ (that is to say in $\Delta'\backslash \Delta$). We can see this in the following way. The first thing to
note if $u \in \dot\Delta$ and $v \in \dot\Delta'$ then neither $uv$ nor $vu$ is in $\Delta \cup \Delta'$. The reason is that 
$2u.v \notin \Z[\tau]$, see Prop.~\ref{isotropy}. So if we assume that $xy \in \dot\Delta'$ then from 
$y = x^{-1}(xy)$ we have that $x^{-1}$ cannot be in $\dot\Delta$. But $x$ and $x^{-1}$ are in $\Delta$, so this shows that $x^{-1}$, and hence also $x$,
is in $K$. Similarly $x = (xy) y^{-1}$ shows that $y$ is in $K$. But then $xy \in K$ since $K$ is a group. This is 
contradicts the assumption that $xy\in \dot\Delta'$, 
so $xy \in \Delta$. This finishes the proof. \qed

\medskip

We note in passing that 
\begin{equation}\label{KDelta}
K \,\dot\Delta \subset \dot\Delta; \; K \,\dot\Delta' \subset \dot\Delta'.
\end{equation}
since $\Delta$ (resp. $\Delta'$) is a group and the products cannot be in $K$, for it is a subgroup. 

\medskip
We now determine basic the structure of $\Sigma$ as seen from its generation by the two groups
$\Delta$ and $\Delta'$. 

The following set $C$ of elements is a set of representatives for the four non-trivial right cosets of $K$ in $\Delta$ :
\begin{equation}\label{leftCosets1}
\left\{\frac{1}{2}(0,1,\tau', \tau), \frac{1}{2}(0,-1,\tau', \tau),\frac{1}{2}(0,1,-\tau', \tau),\frac{1}{2}(0,1,\tau', -\tau)\right\}\,.
\end{equation}
Similarly for  $C'$ and the right cosets of $K$ in $\Delta'$. 

Now, any non-trivial element of $\Sigma$ can be written as an alternating product of elements of $\Delta$
and $\Delta'$. Using the coset representatives we can always write such a product in the form
\begin{equation}\label{standardForm}
w a_m a_{m-1} \dots a_2 a_1 
\end{equation}
where the $a_i$ are alternately in $C$ and $C'$ and $w \in K$, since $K$ is a subgroup of both $\Delta$ and $\Delta'$.

\medskip
We are going to prove that this writing is unique: each element of $\Sigma$ can be written uniquely in this form.
The key to this is to introduce the set $E\subset \Z[\tau]^4/4\Z[\tau]^4$, defined as the set of all 
\begin{equation}\label{setE}
\{ (x_1+y_1\tau, x_2+y_2\tau, x_3+y_3\tau, x_4+y_4\tau) \mod 4\Z[\tau]^4\}
\end{equation}
where the $x_i$ and $y_i$ are integers and
\begin{itemize}
\item all the $x_i$ are congruent to $1$ or $3$ modulo $4$;
\item all the $y_i$ are congruent to $0$ or $2$ modulo $4$;
\item the number of $x_i$ congruent to $1$ is even ;
\item the number of $y_i$ congruent to $2$ is odd.
\end{itemize}

It is straightforward to see that this set has $64$ elements. Using conjugation, we have the second
set $E'$ which is easily seen to have the same four properties except that the number of $x_i$ congruent to $1$
is odd (we are retaining the $x +y\tau$ form). 
These two sets have  remarkable properties: 
\begin{lemma}\label{EEResult}\
\begin{itemize}
\item[\rm (i)] $(1/2) EE =E$ \quad and \quad $(1/2) E'E' =E'$;
\item[\rm (ii)] $E \cap E' = \emptyset$; 
\item[\rm (iii)] $4 C C' \subset E$ \quad and \quad $4C'C \subset E'$.
\end{itemize}
 \end{lemma}
 
Part (ii) is clear from the definitions. It is only necessary to verify the first parts of (i) and (iii).  At the present time we have depended on brute force computation to prove this. For example, in looking at $4CC'$ in (iii), we have the following of
 these $16$ products: 
 \[\begin{array}{cccc}
(1 , -\sqrt{5} , \sqrt{5} , \sqrt{5}) & (3 , -\sqrt{5} , -1 , 1 )&
 (-1 , 3 , \sqrt{5} , -1 )&(-1 , -3 , 1 , \sqrt{5}) \\
(3 , -\sqrt{5} , 1 , -1 )&(
1 , -\sqrt{5} , -\sqrt{5} , -\sqrt{5} )&( 1 , 3 , 1 , \sqrt{5} )&(
1 , -3 , \sqrt{5} , -1) \\
(-1 , -3 , \sqrt{5} , 1 )&(
1 , -3 , -1 , \sqrt{5} )&(
1 , \sqrt{5} , \sqrt{5} , -\sqrt{5} )&( -3 , -\sqrt{5} , 1 , 1) \\
(-1 , 3 , -1 , \sqrt{5} )&(
1 , 3 , \sqrt{5} , 1 )&(
-3 , -\sqrt{5} , -1 , -1 )&(
1 , \sqrt{5} , -\sqrt{5} , \sqrt{5} )
\end{array} \]
 from which one can convert into $x+y\tau$ form to show $4CC' \subset E $. Conjugating
we see that $4C'C \subset E'$.  \qed

\medskip

Using Lemma~\ref{EEResult}, let us show that is not possible for
$z a_m a_{m-1} \dots a_2 a_1$ of \eqref{standardForm} to be equal to $1$, except trivially. 
The $a_i$ alternate between  $C$ and $C'$, and, for definiteness assume that $a_1\in C'$. This analysis of this case
uses $E$, while the case with $a_1\in C$ uses $E'$.

Now suppose that $w a_m a_{m-1} \dots a_2 a_1 =1$. Then $a_m a_{m-1} \dots a_2 a_1 =w^{-1}$.
For this to happen $m\ge 2$.
If $m$ is odd we can take the $a_m$ to the other side to get
$a_{m-1} \dots a_2 a_1 = a_m ^{-1}w^{-1}$. In either case the right hand side is of level at most $1$,
and now we can assume that $m=2n$ is even and $n\ge1$.
 
We  will show that the level of the left-hand side is $n+1$, leading to a contradiction of levels.

Starting with $4a_2 a_1 \in 4CC' \subset E$, which shows that $a_2a_1$ has level $2$. Next we have
\[8a_4a_3a_2a_1 = (1/2) (4a_4 a_3)  (4a_2 a_1) \subset (1/2)EE = E \,,\]
showing that $a_4a_3a_2a_1$ has level $3$. If $n>2$ we can continue in this way to show that
$a_6a_5a_4a_3a_2a_1$ has level $4$, and so on. This shows that indeed the level of the left-hand side
is $n+1$. Of course the argument when $a_1 \in C$ can be obtained by conjugating the equation. 

Now to show the uniqueness of the standard form, suppose that
\[w a_m a_{m-1} \dots a_2 a_1  = z b_p b_{p-1} \dots b_2 b_1 \,\]
are two expressions in the form \eqref{standardForm}.

Supposing that the two expressions are different, we can assume that $p+m$ is minimal for
such an equality. We can push $w$ to the other side, and so assume from the outset that $w=1$.
If $a_1 \in \dot\Delta$ and $b_1 \in \dot\Delta'$, or vice-versa, then we can invert all the elements
of the left hand side and get a new reduced expression that equals $1$, and we know that this cannot happen.

If on the other hand $a_1$ and $b_1$ are both in $\dot\Delta$, say, then 
$b_1a_1^{-1} \in \Delta$ and we can rewrite the right-hand side, to bring it into standard form again, 
while not increasing the number of terms on the right-hand side, and in so doing reduce
the combined length $p+m$, contrary to our minimality assumption. So this case does not happen either.
This proves the uniqueness of expression in \eqref{standardForm} 
Putting this together we have shown:

\begin{tvr}\label{structureOfSigma}
Every element of $\Sigma$ is uniquely expressible in the form \eqref{standardForm}. The group $\Sigma$ is the amalgamation of the two groups $\Delta$ and $\Delta'$ over $K$, i.e. it is the free product of $\Delta$ and $\Delta'$ factored by the normal subgroup generated identifying the elements of $a\in K$ appearing in $\Delta$ and $\Delta'$.
\qed
\end{tvr}

\begin{figure}\centering
\includegraphics[scale=0.35]{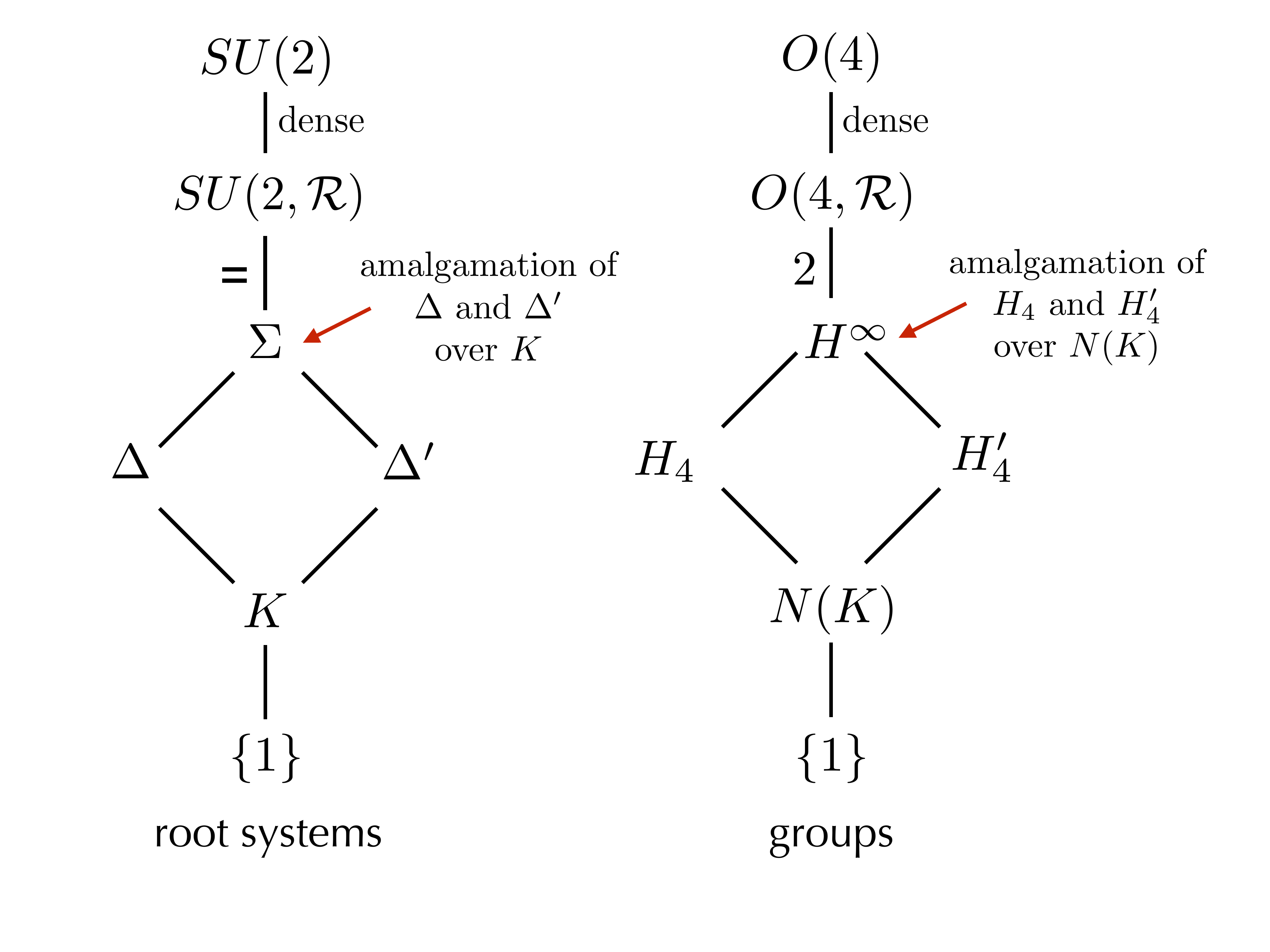}
\caption{How $\Sigma$ and $H^\infty$ appear as amalgamations.}
\label{amalgamation}
\end{figure}
\section{$W(K)$ and its normalizer}

We know that $H^{\infty}[4]$ is generated by two copies of the Coxeter group of type
$H_{4}$, namely $H_4$ and $H_4'$. In this section we shall call them $W(\Delta)$ and $W(\Delta')$.
Both of them contain the subgroup
$W(K)$ generated by the reflections defined by $K$, and this is a Coxeter group (of type $D_{4})$. Using the root bases described above we see that together they create the
Coxeter diagram \eqref{HinftyDiagram}. As we will see below, it appears that there
is an additional relation beyond the Coxeter relations, see \eqref{additionalRelation}. We shall
show the origin of this relation and how together with the obvious Coxeter relations, we
obtain a presentation of $H^{\infty}[4]$.

\medskip

We begin with looking at the group 
\begin{equation}\label{normalizer}
N(K) := \{w \in W(\Delta)\, : \, w(K) = K \} \, .
\end{equation}
Similarly we have $N(K)'$ for the stabilizer of $K$ in $W(\Delta')$.

Since $K$ is a root system of type $D_4$, its complete group of automorphisms is the 
semi-direct product of $W(K)$ and the symmetric group $S_{3}$ (which acts as diagram automorphisms).

\begin{lemma}\label{normalizerK} \
\begin{itemize}
\item[{\rm (i)}]$N(K)$ is the semi-direct product of $W(K)$ and a cyclic group of order $3$.
\item[{\rm (ii)}] $N(K)$ is the  normalizer of $W(K)$ in $W(\Delta)$ . 
\item[{\rm (iii)}] The parallel results hold for $W(\Delta')$ and its normalizer $N(K)'$.
\item[{\rm (iv)}] $N(K) = N(K)'= W(\Delta)\cap W(\Delta') = H_4 \cap H_4'$. 
\end{itemize}
 
\end{lemma}

\noindent
Proof:  The set $K$ is a root system of type $D_4$ and $W(K)$ is its Weyl group.
Thus $N(K)/W(K)$ is a subgroup $S_{3}$ of
diagram automorphisms of $D_{4}$.
Comparing the orders $2^{6}.3^{2}.5^{2}$ of $W(\Delta)$
and $2^{6}.3$ of $W(K)$ (\cite{Bour}), we see that the index $[N(K):W(K)]$ is either $1$ or $3$. 
We shall now see that it is $3$.

The Sylow $3$-groups of $W(\Delta)$ are of order $9$.
They are in fact the direct product of two cyclic groups of order $3$. We can see this from the
fact that there are subgroups of type $A_{2} \times A_{2}$ inside $H_{4}$. E.g. we have
the root pairs  $\{(1,0,0,0), (1/2)(-1,1,1,1)\}$ and $\{u:=(1/2)(0,-1,\tau',\tau), v:=(1/2)(0,\tau',\tau,-1)\}$,
both of which are $A_{2}$ bases. Let $\Gamma$ be the associated $A_{2}\times A_{2}$
root system. 

Multiplying, we find that $uv=( 1/2, -1/2, -1/2, -1/2 ) \in K$ and 
\[r_{u}r_{v}(x_{1},x_{2},x_{3},x_{4})= uv(x_{1},x_{2},x_{3},x_{4})vu = (x_{1},x_{3},x_{4},x_{2})\,.\]
 This three cycle performs a diagram
automorphism of the $D_4$ root system $K$ and together with $W(K)$ it generates $N(K)$.

Although $u,v \in \Delta\backslash K$, their product is in $K$. Thus we have $u',v' \in \Delta'\backslash K$
and their product is the same element of $K$. They serve to produce the very same
three cycle, and hence lie in the normalizer $N(K)'$ of $W(K)$ in $W(\Delta')$.
Thus $N(K) = N(K)'$. 

Any element of $w \in W(\Delta') \cap W(\Delta)$ maps $K= \Delta\cap \Delta'$ into itself
and so stabilizes $K$. Thus it is in $N(K) = N(K)'$.
   \qed
   
\medskip
Let 
\begin{equation}\label{additionalRelation}
s:= r_{u}r_{v}= r_{u'}r_{v'} \,.
\end{equation}
\bigskip

Using \eqref{additionalRelation} we can write down a presentation for the group $H^{\infty}[4]$.
The presentation uses abstract generators $R_{\alpha}, R_{\alpha'}$ and relations
that involve the reflections $r_{\alpha}, r_{\alpha'}$. This is discussed in the next section.

\medskip

\noindent{Comment}: All of the Sylow $3$-subgroups of $W(\Delta)$ are conjugate
by elements of $W(\Delta)$. The ones that are in $N(K)$ are Sylow $3$-subgroups
of it and so are conjugate by elements of $N(K)$ to each other, and in particular
to the one that we used above. These conjugations must conjugate $\Gamma$
into other $A_{2}\times A_{2}$ root systems in $\Delta$. Each of these root
systems generates another relation, just like \eqref{additionalRelation}. Of course
these are just ones that follow from the conjugation of \eqref{additionalRelation}. 

\section{A presentation of $H^{\infty}$}

The purpose of this section is to prove

\begin{tvr}\label{HinfityRelations}
$H^{\infty}= H^\infty[4]$ is generated by the following generators and relations:
\begin{itemize}
\item[{\rm (i)}] generators $R_{\alpha}$, $\alpha \in \Delta$, and 
$R_{\alpha'}$, $\alpha' \in \Delta'$;
\item[{\rm (ii)}] $R_{\alpha}^{2}= 1= R_{\alpha'}^{2}$ for all $\alpha \in \Delta$
and for all $\alpha' \in \Delta'$;
\item[{\rm (iii)}] relations $R_{\alpha}= R_{\alpha'}$ for all $\alpha \in K$;
\item[{\rm (iv)}] relations $R_{\alpha}R_{\beta}R_{\alpha}^{-1} = R_{r_{\alpha}(\beta)}$ for
all $\alpha, \beta \in \Delta$; and similar relations for $\Delta'$;
\item[{\rm (v)}] with $u,v\in \Delta\backslash K$ and  $u',v' \in \Delta'\backslash K$ as given above, $R_{u}R_{v}= R_{u'}R_{v'}$.
\end{itemize}
\end{tvr}

Notice that $R_{\alpha}= R_{\alpha}R_{\alpha}R_{\alpha}^{-1} = R_{-\alpha}$
for all $\alpha \in \Delta$ and similarly for $\Delta'$.

\medskip
The proof depends on working explicitly with the group $H^\infty$ and the reflections $r_\alpha$, 
$r_{\alpha'}$  as they appear in terms of the algebra of  the quaternions, and showing that we can
untangle any relation written in terms of these generators by only using the relations corresponding to 
those itemized in the statement of the Proposition. 
 
We begin by recalling that  for  $a \in \cS$, the quaternionic form of the reflection $r_{a}$ in $a$ \eqref{quaternionReflection}.
  The effect of a product of reflections 
 $ r_{a_{k}}\cdots r_{a_{2}}r_{a_{1}}$ acting on $x$ is

\begin{equation}\label{reflectionProduct}
 x \mapsto \begin{cases}     - a_{k}
  \widetilde{a_{k-1}}\dots \widetilde{a_{2}}a_{1} \,\widetilde x\, a_{1}\widetilde{a_{2}} \dots \widetilde{a_{k-1}}a_{k} \quad \mbox{if $k$ is odd}\\
  \widetilde{a_{k}}
 a_{k-1}\dots a_{2}\widetilde a_{1}\,x\,\widetilde a_{1} a_{2} \dots a_{k-1} \widetilde{a_{k}} \quad \mbox{if $k$ is even}\,.
 \end{cases}
 \end{equation}

When we write out products of reflections as two-sided products, as in \eqref{reflectionProduct},
this involves alternately conjugating elements (this is quaternion conjugation, which
we are designating by tilde), with the explicit form depending on whether there are an even
or an odd number of reflections involved, and a possible overall sign change.  This is rather messy to write down but quite trivial to do in practice. In order to avoid lots of notation we simply signify the whole conjugation process by putting a hat on each symbol, e.g. $\hat a$, which is capable of being conjugated. Thus the symbol beneath the hat
may or may not be conjugated according to the overall length being even or odd. Thus we write

\begin{equation}\label{reflectionProductSimple} r_{a_{k}}\dots r_{a_{1}}(x) = (-1)^{k}\widehat{a_{k}}\dots \widehat{a_{1}} \widehat x \,\widehat{a_{1}}\dots \widehat{a_{k}} \,.
\end{equation}
Notice that conjugation stabilizes each of $K$, $\dot\Delta$, and $\dot\Delta'$.

\medskip

\noindent{Proof:}(Prop.~\ref{HinfityRelations})
Let 
\begin{equation}\label{generalRelation}
r_{a_{k}} \dots r_{a_{2}}r_{a_{1}} = 1 
\end{equation}
be a non-trivial relation in $H^{\infty}$, written in terms of root reflections $r_{a_{j}}$ where the
$a_{j}\in \Delta \cup \Delta'$. 
Our objective is to show that such a relation is a consequence of the relations (i)--(v). Equivalently we wish to show that the word
on the left-hand side of \eqref{generalRelation} can be reduced to the empty 
word by using only these relations. We can rewrite the left-hand side in the form \eqref{reflectionProductSimple}.

To begin with consider any product
$r_{b_{p}}\dots r_{b_{1}}$ with all $b_{j} \in \Delta$ and its 
action on  $\mathbb H$: 
\[r_{b_{p}}\dots r_{b_{1}}(x) = (-1)^{p}\widehat{b_{p}}\dots \widehat{b_{1}} \widehat x \,\widehat{b_{1}}\dots \widehat{b_{p}} \,,\]
for all $x \in \mathbb H$.

Let $\widehat{b_{p}}\dots \widehat{b_{1}} =:B_{L}$ 
and  $\widehat{b_{1}}\dots \widehat{b_{p}} =:B_{R} $. Then $B_{L}, B_{R} \in \Delta$ since 
$\Delta$ is a group, and the mapping assumes the form 
\[x \mapsto (-1)^{p} B_{L}\widehat x\,B_{R}\, .\] 

There seems to be no simple relationship between $B_{L}$ and $B_{R}$, and it seems
quite possible that one would belong to $K$ while the other would not. 
However, if both $B_{L}, B_{R} \in K$ then we can see that 
$r_{b_{p}}\dots r_{b_{1}} \in N(K)$. In fact,  using the fact that the central inversion $-1\in W(K)$
we see that for all $x\in K$, $ r_{b_{p}}\dots r_{b_{1}}(x) = (-1)^{p} B_{L}\widehat x\,B_{R} \in K$,
see \eqref{normalizer}. Notice that the fact that $ r_{b_{p}}\dots r_{b_{1}} \in N(K)$
is a fact that takes place entirely inside the group $W(\Delta)$, that is to say, it is determined
entirely by the structure of $W(\Delta)$ as a Coxeter group. In view
of the structure of $N(K)$, we can rewrite the product $ r_{b_{p}}\dots r_{b_{1}}$ 
so as to assume that all the $b_{j} \in K$ with the possible exception of
$b_{p}, b_{p-1}$ which may be $u,v$ or $v,u$, see \eqref{additionalRelation}. This is all deducible 
from the Coxeter relations of $W(\Delta)$. 

\medskip
We now begin to parse the word on the left-hand side of \eqref{generalRelation}, utilizing the form
of \eqref{reflectionProductSimple}, from right to left. For definiteness
we will assume that $a_{1} \in \Delta$. We now move to the left until we reach 
the first $a_{p+1} \notin \Delta$. If all the roots involved are in $\Delta$, so $p=n$,  there is
nothing to prove since the whole word is in $W(\Delta)$ and the relation \eqref{generalRelation}
is a consequence of the Coxeter relations of $W(\Delta)$. Assuming 
that $p \ne n$, we rewrite $r_{a_{p}}\dots r_{a_{1}}$ using the two blocks
$B^{(1)}_{L} = \widehat{a_{p}}\dots \widehat{a_{1}}$ and 
$B^{(1)}_{R} = \widehat{a_{1}}\dots \widehat{a_{p}}$ which are both in $\Delta$.
In the case that  $r_{a_{p}}\dots r_{a_{1}} \in N(K)$ we rewrite with the word (using the Coxeter
relations of $W(\Delta)$) so that all the $a_j$ involved are in $K$ with the possible exception 
that $a_{p},a_{p-1}$ are the pair $u,v$ in some order. In that case we use \eqref{additionalRelation}
to replace this pair with $u', v'$ in the corresponding order, thus shifting
these last two letters into $\Delta'$ and passing them on to the next block and using only  $\widehat{a_{p-2}}\dots \widehat{a_{1}}$
to make our blocks $B^{(1)}_{L}$ and $B^{(1)}_{R}$. These two blocks
are then elements of $K$. So in the case $r_{a_{p}}\dots r_{a_{1}} \in N(K)\backslash K$, we do this little trick
to make $B^{(1)}_{L}, B^{(1)}_{R} \in K$ and pass along the non-$K$ terms to subsequent blocks.
We write $q=p$ or $q=p-2$ accordingly. 

We may assume that this rewriting and the passing along of pairs $r_{u}r_{v}$ or $r_v r_u$,
has been done in advance. That means that each block multiplies out to be either
in $K$ or in $\dot\Delta$.

\medskip
We now continue the parsing process from where we left off, this time producing a new
set $a_{q}, \dots, a_{p+1}$ which are of elements in $\Delta'$ whose left and right product  
produce blocks $B^{(2)}_{L}$ and $B^{(2)}_{R}$, which are both in $\Delta'$. Again, if these
two blocks are both in $K$ then $r_{a_{s}}\dots r_{ a_{q+1}} \in N(K)' =N(K)$, and this is
derivable from the Coxeter relations of $W(\Delta')$. 

There is a little bit of extra
attention needed here. The elements $a_{j}$ appearing $B^{(2)}_{L}$ and $B^{(2)}_{R}$
have their conjugations determined by the parity of their positions in the sequence of reflections
appearing in \eqref{generalRelation}. But as an element of $W(\Delta')$, the action of $r_{a_{s}}\dots r_{ a_{q+1}}$
might be either $u \mapsto \pm B^{(2)}_{L}\hat u \,B^{(2)}_{R}$, or the conjugations
may be interchanged so its action would be
\[u \mapsto \pm \widetilde{B^{(2)}_{R}}\widetilde{\widehat u}\, \widetilde{B^{(2)}_{L}} \,.\]
In either case, if both  $B^{(2)}_{L}$ and $B^{(2)}_{R}$ are in $K$, the $r_{a_{s}}\dots r_{ a_{q+1}}$
maps $K$ into $K$ and so lies in $N(K)$.

We continue this process until we finally achieve a rewriting of left-hand side of \eqref{generalRelation} in the form
\begin{equation}\label{blockProduct}
(-1)^{n} B^{(h)}_{L} \dots B^{(1)}_{L} (\widehat x)  B^{(1)}_{R}\dots B^{(h)}_{R}  = x\,,
\end{equation}
where $n$ is the original length of the word we began with. This equation is valid for all
$x \in \mathbb H_{4}$. Of these, all blocks multiply to be in $K$ or $\dot\Delta$
or $\dot\Delta'$. 

Now if all the blocks here lie in $K$ then, as we saw above,
the entire effect on $x$ is that of an element of $N(K)$, and the entire relation can be 
determined, block by block, only using the Coxeter relations of $W(\Delta)$ and $W(\Delta')$. 
So let us suppose that at least one block is either in $\dot\Delta$ or $\dot\Delta'$.
Choose a first such block to appear either to the left or right of $x$;  call it $B$. The blocks, if
any, that lie between $x$ and $B$ are in $K$.

For definiteness suppose $B$ is of type $\dot\Delta$ and suppose it is on the right. Choose $x \in \dot\Delta'$.
Then on the left-hand side we have  a partial product in the set
\[ \dots \dot\Delta'  [\mbox{ possible blocks in} \,K]\, \dot\Delta \dots \]
From this we see that every element in the left-hand side has level at least $2$,
see Prop.~\ref{structureOfSigma} and the discussion preceding it. The level can only increase
if there are any other alternations of blocks of types $\dot\Delta$ and $\dot\Delta'$ in \eqref{blockProduct}.
The right hand side of \eqref{blockProduct} has level $1$. This is a contradiction. This argument works in the same way if $B$ is on the left.

This contradiction shows that all the blocks must lie in $K$, and so the 
relations (i)-(v) Prop.~\ref{HinfityRelations} suffice
to reduce the left-hand side of \eqref{generalRelation} to the empty word. This proves that these relations are a presentation of $H^{\infty}$. \qed

\medskip

 There is a more succinct way of saying this
that is useful. Let $\cC$ be a set of coset representatives for all the non-trivial right cosets of $H_4 \mod N(K)$,
and let $\cD$ be a set of coset representatives for all the non-trivial right cosets of $H'_4 \mod N(K)$. 
That is 
\[H_4 = N(K) \cup \bigcup_{y \in \cC}  N(K)\,y\quad \mbox{(disjoint union)} \,\]
and $\cD$ does the same for $H_4'$. A natural choice would be to relate \,$\cC$ and $\cD$ by using
the conjugation $(\phantom{x})'$. However it is better to keep the freedom for other choices, as we shall see below.

\begin{tvr}\label{amalgamation4}\
\begin{itemize}
\item[{\rm (i)}]
Let $\cC$ {\rm (}resp.\,$\cD${\rm )} be a set of coset representatives for all the non-trivial right cosets of $H_4 \mod N(K)$ {\rm (}resp. $H'_4 \mod N(K)${\rm )}. Then very element of $H^\infty[4]$ is uniquely expressible in the form
	\begin{equation}\label{ue}
		z\, y_k \dots y_1
			\end{equation}
where $z \in N(K)$, and where the $y_j$ alternate between $\cC$ and $\cD$. 
\item[{\rm (ii)}] $H^{\infty}$ is the amalgamation of $H_{4}=W(\Delta)$ and $H_{4}'=W(\Delta')$ 
over their intersection $N(K)$. 
\end{itemize}
 \qed
\end{tvr}

\begin{tvr}\label{amalgamation3}\
Let us choose $\cC$ and $\cD$ of Prop.~\ref{amalgamation4} so that the representatives are actually in $H_3$ (resp. $H'_3$) when the coset contains elements of these subgroups. Then
\begin{itemize}
\item[{\rm (i)}]
every element of $H^\infty[3] \subset H^\infty[4]$ is uniquely expressible in the form
	\[ z\, y_k \dots y_1	\]
where $z \in H_3\cap H_3'$, and where the $y_j$ alternate between $\cC$ and $\cD$\,; 
\item[{\rm (ii)}] $H^{\infty}[3]$ is the amalgamation of $H_{3}=W(\Delta[3])$ and $H_{3}'=W(\Delta[3]')$ 
over their intersection. 
\end{itemize}
\end{tvr}
\noindent
Proof: Notice that $K \cap V = \{(0,\pm 1 , 0,0), (0,0, \pm 1, 0),(0, 0,0,\pm 1 )\} = \Delta[3] \cap \Delta[3]'$, which is a root
system of type $A_1\times A_1 \times A_1$ inside $\Delta[3]$. The stabilizer of $K \cap V$ in $H_3$
is the semi-direct product of the Weyl group of $A_1\times A_1 \times A_1$ and the cyclic group of
order three arising by cycling the last three components of $K\cap V$ (or equivalently cycling the three $A_1$s).  
All of these elements are in $N(K)$ and from this we see that the stabilizer
of $K \cap V$ in $H_3$ is $N(K) \cap H_3$.  In the same way the stabilizer of 
$K \cap V$ in $H'_3$ is $N(K) \cap H'_3$, and from this we see 
\[N(K) \cap H_3 = N(K) \cap H'_3 = H_3 \cap H'_3 \,.\]

Now $H^\infty[3]$ is the subgroup of $H^\infty[4]$ generated by $H_3$ and $H'_3$.
According to Prop.~\ref{amalgamation4} every element of $H^\infty[4]$ is uniquely
in the form \eqref{ue}.
Amongst all the expressions \eqref{ue}, consider all those in which 
all the elements occurring are actually in $H_3$ or $H'_3$, and $z\in H_3 \cap H'_3 \subset N(K)$.
All of these expressions are obviously elements of $H^\infty[3]$. What's more, every element of
$H^\infty[3]$ can evidently be written in this form. So $H^\infty[3]$ is composed precisely
of all of these products. Since the expression in this form is unique, we see that in fact it says that 
this is the amalgamation of $H_3$ and $H_3'$ over their intersection. \qed

\section{The orthogonal group}\label{orthoGroup}

Define $O(\Sigma)= O(\Sigma[4])$ to be the group of all orthogonal transformations of $\mathbb R^{4}$
that stabilize the root system $\Sigma$. We already know a lot of transformations
in $O(\Sigma)$, namely all the elements of $H^{\infty}$. Along with $O(\Sigma)$
we shall consider the orthogonal group $O(4,\cR)$ of all orthogonal 
matrices with coefficients in $\mathcal R$. Since $\Sigma \subset \mathcal R^{4}$
and it contains a standard orthogonal basis of $\mathcal R^{4}$, we have
$O(\Sigma) \subset O(4,\cR)$.

Notice that there is an orthogonal transformation of $\Sigma$ that interchanges
$\Delta$ and $\Delta'$, namely the automorphism of $\Sigma$ that interchanges the
last two coordinates $(x_{1},x_{2},x_{3},x_{4}) \mapsto (x_{1},x_{2},x_{4},x_{3})$. Call this
transformation $\iota$. (Of course it has numerous conjugates.) Evidently
$O(\Sigma)$ contains the subgroup $\langle H^{\infty}[4], \iota\rangle$ generated
by $H^{\infty}$ and $\iota$. Now we prove the reverse inclusion.

\begin{tvr}\label{orthoGroupRevealed} \
\begin{itemize}
\item[\rm{(i)}] $O(\Sigma[4]) = O(4,\cR) $ \,;
\item[\rm{(ii)}] $O(\Sigma[4]) = \langle H^{\infty}[4], \iota\rangle = \langle \iota \rangle \ltimes H^\infty[4] =
\langle H_4, \iota\rangle$\,;
\item[\rm{(ii)}] $[O(4,\cR)):H^{\infty}[4]]=2 \,.$
\end{itemize}
\end{tvr}

\noindent Proof: (i) Let $g\in O(4,\cR)$. Then $g$ stabilizes $\mathcal S \cap \mathcal R^{4}
= \Sigma$, and so $g \in O(\Sigma)$.

(ii) We know that $K$ is a root system of type $D_{4}$, for which
\eqref{D4Base} is a root base. Let $B$ be any $D_{4}$ root base in $\Sigma$.

We will try to use elements of $H^{\infty}$ to bring $B$ back to the 
basis  \eqref{D4Base}. We will see that we can almost do this. To finish the job
we may need to reverse the last two coordinates, hence use the transformation $\iota$.

Let $V:=(0,\mathbb R, \mathbb R, \mathbb R)$, as before, and recall
that $V\cap \Sigma[4] = \Sigma[3]$, see \eqref{Sigma3FromSigma4}. 

This base has the traditional $D_{4}$ Coxeter diagram with one central node
and three nodes attached to it.
Using $H^{\infty}$ and simple transitivity, we can assume that one of the non-central
nodes of the base
elements of $B$ is $(1,0,0,0)$. Then the other two non-central nodes are in $V\cap \Sigma[4]$
and so in $\Sigma[3]$. We take one of these and again using transitivity--- this time 
of $H^{\infty}[3]$--- to bring this node to $(0,1,0,0)$. The remaining non-central node
is now of the form $(0,0,*,*)$. The only elements of $\Sigma[3]$ of this form are
$(0,0,\pm 1, 0)$ and $(0,0,0,\pm 1)$. We can use $\iota$ to assume that it is
$(0,0,\pm 1, 0)$. Since the sign changes are in $H^{\infty}[3]$, we can assume
this element is $(0,0,1,0)$. There remains only the central node. Its scalar product with 
the three other nodes (that is, with $(1,0,0,0), (0,1,0,0), (0,0,1,0)$ is $-1/2$ and so it has
the form $(1/2)(-1,-1,-1,*)$. Given that it is a vector of length $1$, the last term must
be $\pm 1$. We again have the choice of sign.

The upshot of this is that we have brought the base $B$ to the standard basis of
\eqref{D4Base} using $H^{\infty}$ and $\iota$. Since any orthogonal transformation of
$\mathbb R^{4}$ is determined entirely by its action on any base, we have proved
that $O(\Sigma) \subset \langle H^{\infty}, \iota\rangle$. The rest of (ii) follows once we know
that (iii) is true.

(iii) We know that $\iota^{2}=1$ and $\iota$ 
stabilizes $\Sigma$ and hence normalizes $H^{\infty}$ (which is generated
by the reflections in the elements of $\Sigma$). Also $\iota$ interchanges $\dot\Delta$
and $\dot\Delta'$ (since it is an odd permutation of the coordinates) while stabilizing $K$.
It follows that conjugation by $\iota$ interchanges $W(\Delta)= H_4$ and $W(\Delta')= H_4'$ while stabilizing
their intersection, which is $N(K)$. This means that in choosing the coset representatives
used in Prop.~\ref{amalgamation4} we can choose $\cD = \iota(\cC)\iota^{-1}$. We will assume this. 

Now we prove that $\iota \notin H^\infty$. If, on the contrary, it were in $H^\infty$ then using
Prop.~\ref{amalgamation4} we could write $\iota$ uniquely in the form
$\iota = z y_k \dots y_1 $
where $ z\in N(K)$
and the other elements are alternately in $\cC$ and $\cD$. 
Then
\begin{align*}
\iota &= \iota \iota \iota^{-1} = \iota z y_k \dots y_1 \iota^{-1}\\
&=  \iota z \iota^{-1} \iota y_k\iota^{-1} \dots \iota y_1 \iota^{-1}\,.
\end{align*}
This gives us a second way of writing $\iota$ as an element of $H^\infty$, in the form
described in Prop.~\ref{amalgamation4} and accordingly 
they are identical. This obviously is not true if any of the coset terms are present, since $\cC$ and $\cD$ are reversed
by the conjugation and $y_1 \ne \iota y_1 \iota^{-1}$, etc. Thus we end up with 
$\iota = \iota z \iota^{-1}$ and so $\iota = z \in N(K)$.

Now certainly $\iota \notin W(K)$, for it effects a diagram automorphism of $D_4$, and we have already noted that 
$[N(K):W(K)]=3$. Thus $\iota \notin N(K)$, and this contradiction shows that $\iota \notin H^\infty$.
\qed

\medskip
There is a corresponding result for the three-dimensional case involving $\Sigma[3]$ and $O(3, \cR)$.
Here we take $O(3, \cR)$ to be the subgroup of $O(4, \cR)$ which fixes the vector $(1,0,0,0)$ in $\HH$,
and define $O(\Sigma[3])$ to be the subgroup of $O(3, \cR)$ stabilizing $\Sigma[3]$.
\begin{tvr}\
\begin{itemize}
\item[\rm{(i)}] $O(\Sigma[3]) = O(3,\cR) $ \,;
\item[\rm{(ii)}] $O(\Sigma[3]) = \langle H^{\infty}[3], \iota\rangle = \langle \iota \rangle \ltimes H^\infty[3] =
\langle H_3, \iota\rangle$\,;
\item[\rm{(ii)}] $[O(3,\cR)):H^{\infty}[3]]=2 \,.$
\end{itemize}
\end{tvr}

\noindent Proof: The argument is essentially a repeat of what we saw in Prop.~\ref{orthoGroupRevealed}. 
Since $\Sigma[3] = \cS[3] \cap \cR^3$, it is clear that $O(3,\cR)$ stabilizes $\Sigma[3]$, and this proves
(i). We need only prove that $O(3,\cR) \subset \langle H^\infty[3], \iota\rangle$.
 Working in $\HH$,
$\Sigma[3]$ contains the standard basis vectors $v_1:=(0,1,0,0),v_2:=(0,0,1,0), v_3:=(0,0,0,1)$. Let $g\in O(3,\cR)$
and let $g(v_1) = x \in \Sigma[3]$. Since $H^\infty[3]$ is transitive on $\Sigma[3]$, there is
an element $h \in H^\infty[3]$ so that $hg(v_1) = v_1$. Then $hg$ stabilizes the span 
of $v_2$ and $v_3$ as well as mapping them into elements of $\Sigma[3]$. However, the only
elements of $\cR^3$ of this form are $(0, 0,\pm 1, \pm1)$. Since the sign changes on individual coordinates
are elements of $H^\infty[3]$, we can assume that we have the positive signs here. Thus 
$hg$ fixes both $v_2$ and $v_3$, in which case $hg = 1$, or $hg$ interchanges $v_2$ and $v_3$, in which
case $hg = \iota[3]$.  \qed

\section{Approximation in $SU(2)$}

The unitary group $SU(2)$ plays a significant role in the theory of quantum computation as the group of 
admissible transformations of a quantum qbit. In this section we offer an algorithm by which any element of $SU(2)$ can be 
approximated to any degree of accuracy by the repeated use of just a few very simple fixed elements in $SU(2,\cR)$. In fact, this number can be reduced to five, these five being related to a root base of $\Sigma[4]$.
 Passing to $SO(3)$ via \eqref{SU2toSO3} we can similarly approximate any rotation to any degree of accuracy by elements of $SO(3,\cR)$. 
 
 The idea is this. Inside the quaternion algebra $\HH$ the unit sphere $\cS$ is a group under
 quaternion multiplication and this group is $SU(2)$ from the point of view of the matrix representation of $\HH$ 
 of \eqref{matrixRepOfQuaternions}.
 The full set of roots $\Sigma[4]$ that we have generated out of $\Delta[4]$ and $\Delta'[4]$
 is in fact the subgroup of $\cS$ generated by them and it is the orbit of $\Delta[4] \cup \Delta'[4]$ under $H^{\infty}[4]$.
 Furthermore it is dense on the sphere $\cS$, Prop.~\ref{denseness}. Hence appropriate sets of mirrors (reflecting hyperplanes) associated with the reflections in these elements
partition $\cS$ into convex regions of arbitrarily small diameter. We call such regions {\em chambers}. 

A natural choice is a set \,$\cU_{n}:= \bigcup_{k\le n} \dot\cS_{k} \cup \dot\cS_{k}'$ and we can use the associated hyperplanes
of these points to define the set of chambers that we use. So we note that the chambers here are not absolute things, but rather the outcome of some selection of a finite number of root hyperplanes.
Assuming that the chambers are selected and that one point $p$ has been chosen in one of these chambers, say  $p\in C$, then the idea is that any element 
$x \in\cS$  can be reflected into $C$ using only reflections in these mirrors, 
whereby it is a close approximation of $p$. Then reversing the operations we can take $p$ to some point
equally close  to $x$. 

So we take some set $\cU_{n}$ and use its points to define a set of chambers. We should note that the situation with these chambers is not the same as in the customary case of root systems. The chambers are not all isometric copies
of each other, nor are they necessarily simplices. The cellular decomposition comes about through the sets of reflecting hyperplanes that we use, but we never use all of them, only a finite number of them, and these finite sets are not invariants of $H^{\infty}$.\footnote{If the process is carried out to the end, allowing $n\to \infty$, then every point $y$ not on one of the reflecting hyperplanes is
uniquely specified by the signs that it makes with each of the roots: $y.a$ is either positive or negative. In this
way it looks similar to a Dedekind cut where each real number is specified by the rationals that are less
than it or greater than or equal to it.} 

For $a \in \cS$, the reflection $r_{a}$ in $a$ has the quaternionic form
\[ r_{a}x = x - 2(x.a)a = x - (x\widetilde a + a\widetilde x) a = - a\widetilde x a \, \,\]
 see \eqref{quaternionReflection}.
 Suppose we wish to approximate $x \in \cS$ using the reflections determined by the elements of $\cU_{n}$. The idea is to reflect $x$ a number of times so that it ends up as close as possible to $1 \in \mathbb H$. Choose a chamber
$C_{0}$ containing $1$ (actually $1$ will be at the vertex of such a chamber) and choose a point $z$ close to $1$
in the interior of $C_{0}$. For each element pair $\pm a \in \cU_{n}$ we choose the one which satisfies $a.z < 0$
and let $\cU^{-}_{n}$ denote this subset chosen in this way.  So $C_{0} = \{u \in \cS: u.a\le 0\; \mbox{for all} \;a \in \cU^{-}_{n}\}$.

Now if there is a reflection hyperplane for which
$x$ and $z$ lie on opposite sides, then $a.x >0$ and we see that
\[(r_{a}x-z).(r_{a}x-z) = (x-z).(x-z) + 4 (a.x)(a.z) < (x-z).(x-z) \,,\]
thus moving the point $x$ closer to $z$. We can always do this so that we choose $a$ for which
$a.x$ is as large as possible. 

Almost surely, this process has to stop in a finite number of steps with the point ending up in $C_{0}$. Here is why. 
Suppose that $x = x_{0},x_{1}, x_{2}, \dots $ were an infinite sequence of successive steps applied to the starting point $x$. It has to have limit points, and we choose
one, say $y$,  that is closest to $z$. Now $a.z \le 0$ for all $a\in \cU^{-}_{n}$, and it follows that $y\in C_{0}$.
In fact every iteration brings us closer to $z$ and if $y\notin C_{0}$ there is an $a\in \cU^{-}_{n}$ with
$a.y>0$ and $r_{a}(y)$ is closer to $z$, a contradiction. Almost surely (in the sense of Lebesgue measure) $x$ is not on the boundary of a mirror
so neither is $y$. Thus $y$ is in the interior of $C_{0}$. This means that almost surely there is a $k$ for which 
$x_{k}$ is in the interior of $C_{0}$, and the iteration will stop as soon as it reaches such a $k$. 

From the point of view of approximation theory, we can always assume that $x$ is not on a boundary.
At the end of the process we will have achieved a relation of the form
\[ r_{a_{k}}\cdots r_{a_{2}}r_{a_{1}}(x) = p \in C_{0} \,.\]
If $d$ is the diameter of $C_{0}$, then since each reflection is an isometry of $\mathbb H$ relative to its
usual norm, 
\[ |x- r_{a_{1}}r_{a_{2}}\cdots r_{a_{k}}1| = | r_{a_{k}}\cdots r_{a_{2}}r_{a_{1}}(x) -1| \le d \,.\]
This then gives us the desired approximation
\begin{equation}\label{matrixProduct}
 x \simeq \begin{cases}     - a_{1}
  \widetilde{a_{2}}\cdots \widetilde{a_{k-1}}a_{k}a_{k}\widetilde{a_{k-1}} \cdots \widetilde{a_{2}}a_{1} \quad \mbox{if $k$ is odd}\\
  \widetilde{a_{1}}
 a_{2}\cdots a_{k-1}\widetilde a_{k}\widetilde a_{k} a_{k-1} \cdots a_{2} \widetilde{a_{1}} \quad \mbox{if $k$ is even}\,.
 \end{cases}
 \end{equation}
 
 Treated as elements of $SU(2)$, the matrix entries all come from the ring of dyadic golden numbers $\cR$.
 
To write this fully in terms of the basic reflections in the roots of $\dot \Delta \cup K\cup\dot\Delta'$ we need
 to write each of the roots appearing in $\cU_{n}^-$ in terms of these elements. Thus if $a\in U_{n}^-$
 is $a= r_{b_{m}}r_{b_{m-1}}\cdots r_{b_{1}}b_{0}$ with all the $b_{j}\in \dot \Delta \cup K\cup\dot\Delta'$, then
 \[a = \begin{cases}    \widetilde{b_{m}}
 b_{m-1}\cdots b_{1}\widetilde b_{0} b_{1} b_{2} \cdots b_{m-1} \widetilde{b_{m}} \quad \mbox{if $m$ is odd} \\
 
  - b_{m}  \widetilde{b_{m-1}}\cdots \widetilde{b_{1}} b_{0}\widetilde{b_{1}} \cdots \widetilde{b_{m-1}}b_{m} \quad \mbox{if $m$ is even}
 \,.
 \end{cases}\]

\medskip
 
\noindent {\bf Example:} Here we give some sample results from a computer implementation of the approximation process.
The set $\cU_{3}$ of roots formed from the sets $\dot \cS_{k},\dot\cS_{k}'$ is evaluated for $k=1,2,3$, resulting in a total of 
\[52416 = 2.3.(2^{5} + 2^{9} + 2^{13})\] elements.
 A random element of $SU(2)$ is computed, e.g.
 \[ u = \begin{bmatrix}
 -0.244828 - 0.155561 i &
 0.731977 - 0.616498 i\\
 -0.731977 - 0.616498 i &  -0.244828 +  0.155561 i
 \end{bmatrix} \]
 and the process then determines (in this case six) reflections from roots of $\cR_{3}$ 
 that bring this point on $\cS[4]$ as close to the identity element $1 \in \mathbb H$ as possible (using reflections from
 $\cR_{3}$). These reflections
 are rewritten as products of reflections from the original set $\dot\Delta \cup \dot\Delta'$, giving rise to 
 a total of fifteen.  The matrix 
\[\frac{1}{2}\begin{bmatrix}
  1 + \tau i &  (1 - \tau) i \\
   (1 - \tau) i &1 - \tau i
 \end{bmatrix} \]
 is a typical looking example.
 
After computing the resulting matrix product according to \eqref{matrixProduct}, we obtain the approximation
\[\begin{bmatrix}
-0.245712 - 0.15869 i &
 0.730844 - 0.616694 i \\
   -0.730844 - 0.616694i & -0.245712 + 
  0.15869 i 
  \end{bmatrix} \, ,\]
  the difference being
  \[\begin{bmatrix}
0.000883332 + 0.00312853 i &  0.00113294 + 0.000195296 i\\
 -0.00113294 + 0.000195296 i & 0.000883332 - 0.00312853 i \,
  \end{bmatrix}.\]
  
  This example is entirely typical both in the number of reflections required and in the 
  degree of approximation. 
  If we wished we could rewrite the reflections in the roots of $\dot\Delta$ and $\dot\Delta'$
  in terms of the five generators of \eqref{Coxeter generators}. Computing out more roots, e.g. $\cU_{4}$
  will result in more accuracy. We do not have any decent estimates on the diameters of the chambers
  arising from $\cU_n$. 
  
 \section*{Acknowledgements}
 
 We thank the referee for his/her careful parsing of the manuscript and helpful historical suggestions. 
 
 RVM was supported by Discovery Grant
 461292003, Natural Sciences and Engineering  Research Council of Canada.
 
 \noindent
 JM was supported by Grants-in-Aid for Scientific Research 26400005, 17K05158
(Monkasho Kakenhi, Japan).
 

\end{document}